\documentclass[11pt]{article}

\usepackage{amssymb}
\usepackage{amsthm}
\usepackage{amsmath}
\usepackage{graphicx}
\usepackage{psfrag}
\usepackage{rawfonts,color}

\newtheorem{theorem}{Theorem}
\newtheorem{lemma}{Lemma}
\newtheorem{proposition}{Proposition}
\newtheorem{corollary}{Corollary}

\theoremstyle{definition}
\newtheorem{definition}{Definition}
\newtheorem{remark}{Remark}

\newcommand\al{\alpha}
\newcommand\eps{\varepsilon}
\newcommand\ka{\kappa}
\newcommand\la{\lambda}
\newcommand\om{\omega}
\newcommand\Om{\Omega}
\newcommand\ph{\varphi}
\newcommand\si{\sigma}

\newcommand\R{\mathbb{R}}
\newcommand\C{\mathbb{C}}

\newcommand\N{\mathbb{N}}
\newcommand\RnR{\R^{n} \times \R}
\newcommand\RnI{\R^{n} \times (0,\infty)}
\newcommand\RtR{\R \times \R}
\newcommand\RtI{\R \times (0, \infty)}

\newcommand\Czz{C^{oo}}
\newcommand\Czzsym{C^{oo}_{\rm sym}}
\newcommand\Ssym{S_{\rm sym}}
\newcommand\Spsym{S'_{\rm sym}}

\newcommand{\req}[1]{(\ref{eq:#1})}

\newcommand{\f}{\mathbf}
\newcommand{\fs}{\boldsymbol}

\newcommand{\Crec}{\mathcal{C}_{\rm rec}}
\newcommand{\Drec}{D_{\rm rec}}
\newcommand{\Vrec}{V_{\rm rec}}

\DeclareMathOperator{\Io}{\cal I}
\DeclareMathOperator{\Jo}{\cal J}
\DeclareMathOperator{\Ko}{\cal K}
\DeclareMathOperator{\Mo}{\cal M}

\DeclareMathOperator{\Po}{\cal P}
\DeclareMathOperator{\Qo}{\cal Q}
\DeclareMathOperator{\Ro}{\cal R}

\DeclareMathOperator{\inj}{\rm i}
\DeclareMathOperator{\injS}{\rm i_{S}}
\DeclareMathOperator{\injL}{\rm i_{L^1}}
\DeclareMathOperator{\injC}{\rm i_{\Czz}}

\DeclareMathOperator{\ft}{ \rm F }
\DeclareMathOperator{\ift}{ {\rm F}^{-1} }
\DeclareMathOperator{\cost}{ \rm C }
\DeclareMathOperator{\icost}{ {\rm C}^{-1} }
\DeclareMathOperator{\ftp}{ {\rm F}^\prime }
\DeclareMathOperator{\iftp}{ ({\rm F}^\prime)^{-1} }
\DeclareMathOperator{\costp}{ {\rm C}^\prime }
\DeclareMathOperator{\icostp}{ ({\rm C}^\prime)^{-1} }

\newcommand\norm[1]{\|#1\|}
\newcommand\abs[1]{|#1|}
\newcommand\set[1]{\{#1\}}
\newcommand{\ip}[2]{\langle #1,#2 \rangle}
\newcommand\tm{\subset}
\newcommand{\lap}{\Delta}

\definecolor{farba}{rgb}{0.3,0.1,0.1}
\definecolor{farbb}{rgb}{0.1,0.15,0.2}
\definecolor{farbc}{rgb}{0.1,0.2,0.1}

\usepackage[dvips,
  pdftitle = {Thermoacoustic CT},
  pdfauthor = {Markus Haltmeier und Fidler Thomas},
  bookmarksopen,
  colorlinks,
  linkcolor = black,
  urlcolor  = blue,
  citecolor = black,
  menucolor = farba
]{hyperref}

\begin{document}

\title{Mathematical Challenges Arising in Thermoacoustic Tomography with
Line Detectors}

\author{
Markus Haltmeier\thanks{
\href{mailto:markus.haltmeier@uibk.ac.at}{\tt
Markus.Haltmeier@uibk.ac.at}.} \and Thomas Fidler\thanks{
\href{mailto:Thomas.Fidler@uibk.ac.at}{\tt
Thomas.Fidler@uibk.ac.at}} }

\date{Department of Mathematics, University of Innsbruck,
        Technikerstrasse 21a, A-6020 Innsbruck, Austria,
        \\[1ex]
        Revision, April 2008 (original version: August 2006)
}

\maketitle

\begin{abstract}
Thermoacoustic computed tomography (thermoacoustic CT) has the
potential  to become a mayor non-invasive  medical imaging method.
In this paper we derive a general mathematical framework of a novel
measuring setup introduced in [P. Burgholzer, C. Hofer, G. Paltauf, M. Haltmeier, and
O. Scherzer, \emph{Thermoacoustic tomography with integrating area and line detectors},
IEEE Transactions on Ultrasonics, Ferroelectrics, and Frequency Control, 52 (2005)],
that uses line shaped detectors instead of the usual point like ones. We show that
the {\em three dimensional} thermoacoustic imaging problem reduces
to the mathematical problem of reconstructing the initial data of
the {\em two dimensional} wave equation from boundary measurements
of its solution. We derive and analyze an analytic reconstruction
formula which allows for fast numerical implementation.
\end{abstract}

\noindent {\small {\bf Keywords.} Thermoacoustic; Tomography; Line
detectors; Wave equation; Optoacoustic; Photoacoustic; Limited data;
Circular Radon transform.}
\medskip

\noindent {\small {\bf AMS Classification.} 35L05, 44A12, 65R32.}

\section{Introduction and Motivation}\label{sec:intro}

Thermoacoustic CT  (also called opto- or photoacoustic CT) is a
novel hybrid non-invasive imaging method with applications in
different areas, e.g. in medical diagnostics \cite{KruEtAl03} or
imaging of small animals \cite{WPKXSW03}. It is based on the
excitation of acoustic waves inside an investigated object when
exposed to non-ionizing electromagnetic radiation \cite{HSS05,
GusKar96} and combines the advantages of purely optical imaging
(high contrast) with ultrasonic high resolution \cite{XuMFenWan03,
AKO01}.

So far the generated thermoacoustic waves are measured with several
ultrasonic transducers located outside the illuminated sample. The
measured output of the ultrasonic transducers has been identified
with the restriction of the thermoacoustic pressure field to a
surface enclosing the object. Based on this point-data approximation
the absorption density function can be reconstructed by solving the
problem of recovering a function from its mean values over spherical
surfaces \cite{FinPatRak04, XuMFenWan03}. When conventional
piezoelectric ultrasonic transducers \cite{XuMFenWan03, KruEtAl03}
are used to approximate point-data, the necessity of using small
detectors with high bandwidth is technically hardly realizable
\cite{CoxEtAl04, PalSchGus96}.

To obtain high resolution the size of the detectors has to be taken
into account when modeling the corresponding forward operator. Exact
reconstruction formulas incorporating the detector size have been
derived for large planar detectors in spherical geometry
\cite{art04:HalEtAl} and line detectors with cylindrical circular
recording geometry \cite{BurEtAl05, PNHB07}.

In this paper we establish a mathematical foundation of thermoacoustic CT using
line detectors in general recording geometry. We show that the {\em three dimensional} absorption density function can be reconstructed
by solving the mathematical problem of recovering the initial data
of the {\em two dimensional} wave equation from its solution,
measured by an array of line  detectors, on a recording curve $\Crec$.
Therefore, line detectors offer a reduced numerical  complexity.
Finally three dimensional image reconstruction is achieved by rotating
the array  around a single axis.

In the case where the recording curve $\Crec$ is a line, the initial
data of the two dimensional wave equation can be recovered by a {\em
Fourier reconstruction formula} \cite{KoeBea03, KoeEtAl01, PNHB07}.
In this article we present a rigorous mathematical analysis of this
formula, taking into account that the restriction of the solution of
the wave equation to a line is not necessarily absolutely
integrable.

The paper is organized as follows. In Section \ref{sec:tct} we
recall the basic formulas of thermoacoustic CT and introduce the concept of line
detectors. Section \ref{sec:rot} is devoted to the decomposition of
the three dimensional forward operator into a system of two
dimensional operators corresponding to the two dimensional wave
equation. In Section \ref{sec:linear} we deal the case of linear recording
geometry and present the analysis of the Fourier reconstruction formula.

\section{Mathematical modeling of thermoacoustic CT with line detectors}\label{sec:tct}

The basic principle of thermoacoustic CT is the generation of acoustic waves
within an object by illuminating it with non-ionizing pulsed electromagnetic energy.

Assume that a short electromagnetic pulse is emitted
into a weakly absorbing medium at time $\hat{t} = 0$. The absorbed
energy induces thermal heating, which causes thermoelastic expansion
and thereby an initial pressure distribution. This actuates small vibrations
in the medium, which results in sound waves.

The induced pressure field is mathematically modeled by a function
$p: \R^3 \times  [0, \infty) \to \R$, where $p(\f x, \hat{t})$
represents the acoustic pressure at position $\f x \in \R^3$ and
time $\hat{t} =: t/c  \geq 0$. Here $c$ denotes the {\em speed of
sound} which is assumed to be constant. In this case the time
varying pressure field satisfies the three dimensional wave equation
(see \cite{GusKar96, XuMFenWan03, HSS05})
\begin{equation}\label{eq:wave3d}
    \left( \frac{\partial^2}{\partial t^2} - \lap_{\f x} \right) p(\f x,t) = 0\,,
    \quad(\f x,t) \in \R^3 \times(0, \infty)
\end{equation}
with initial conditions
\begin{eqnarray}
    p(\f x,0)  =  f(\f x) \,, & \quad \f x \in \R^3  \,,&    \label{eq:ini3d1} \\
    \frac{\partial p}{\partial t}( \f x, 0)  =  0 \,,  & \quad \f x  \in \R^3 \;.&
\label{eq:ini3d2}
\end{eqnarray}
The mathematical task in thermoacoustic CT is to reconstruct the
initial data $f$ using data of the solution of the three dimensional
wave equation (\ref{eq:wave3d}), (\ref{eq:ini3d1}),
(\ref{eq:ini3d2}) gathered with several detectors located outside
the investigated sample.
\begin{psfrags}
    \psfrag{Receiver}{transducer}
    \psfrag{Absorbers}{absorbers}
    \psfrag{Heat}{Heat}
    \psfrag{deposited}{deposited}
    \psfrag{Illumination}{Illumination}
    \begin{figure}
        \begin{center}
            \includegraphics[width = 0.3\textwidth ]{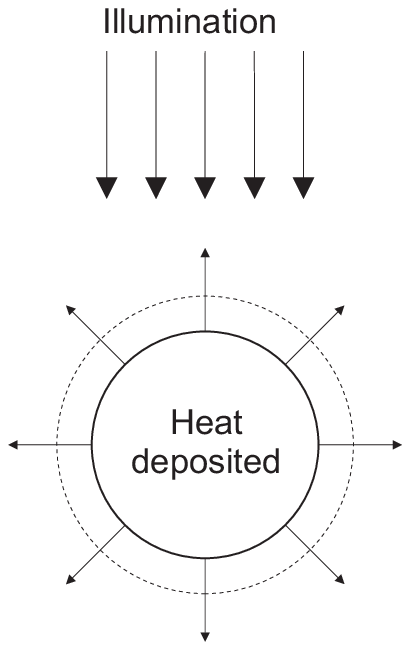} \hfill
            \includegraphics[width = 0.5\textwidth ]{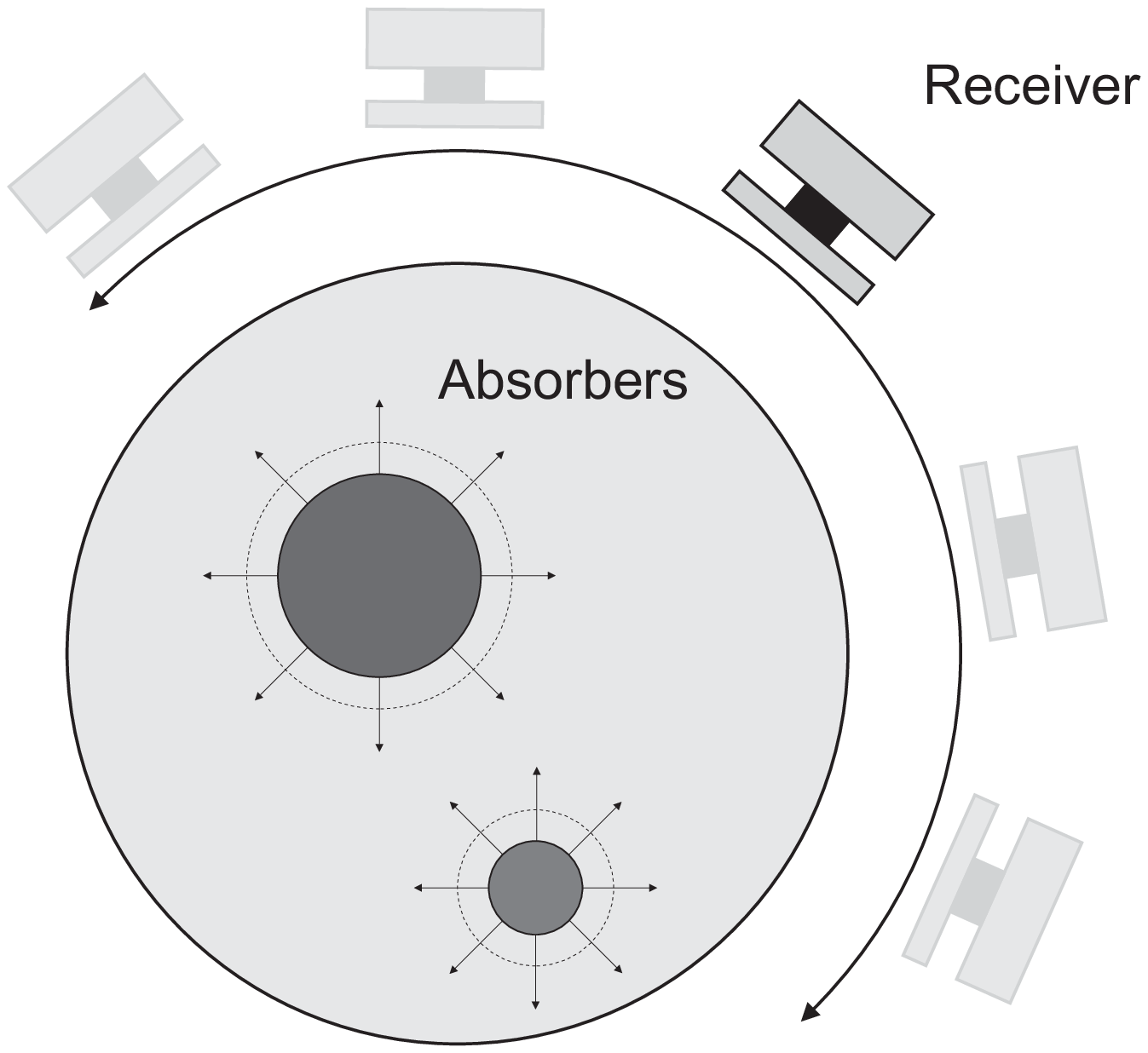}
        \end{center}
        \caption{{\bf Conventional thermoacoustic scanning system.}
                 The non-uniform absorbtion causes outgoing acoustic waves
                 that are measured with an ultrasonic transducer
                 located at different positions.}
        \label{fg:setup}
    \end{figure}
\end{psfrags}

The measured output of a conventional ultrasonic transducer is the
total pressure
\begin{equation*}
       \frac{1}{ \abs{ S }}\int_{S} p(\f x,t)\, d S(\f x)
\end{equation*}
acting on the surface $S$ of the ultrasonic transducer convolved in
time with the temporal impulse response function of the detector
\cite{Jen99}. Here $d S$ denotes the two dimensional surface
measure and  $\abs{ S }$ the surface area of $S$.
In thermoacoustic CT the temporal impulse response function of the ultrasonic
detector is assumed to be an approximation of the
$\delta$-distribution and the surface area $\abs{S}$ to be small. In
this case the measurement data  can be approximated by
\begin{equation*}
        p(\f x_{\rm rec}, t) \cong  \frac{1}{ \abs{S}} \int_{S} p(\f x,t)\, d S(\f x)
\end{equation*}
where $\f x_{\rm rec}$ is the center of the surface $S$.

In practice conventional piezoelectric ultrasonic transducers placed
on a recording surface $S_{\rm rec}$ are used to approximate
point-like detectors \cite{KruEtAl03, XuMFenWan03}. In this case
thermoacoustic CT leads to the mathematical problem of recovering a
function from its integrals over spheres with centers $\f x_{\rm
rec}$ (locations of the {\em point-like} detectors) on $S_{\rm rec}$
\cite{XuMFenWan03, FinPatRak04, KoeEtAl01}. Every ultrasonic
transducers has a finite size (typically in the range of mm) and
therefore algorithms that are based on the assumption that point
measurement data is collected give blurred reconstructions.

{\em Line detectors} measure the total pressure
\begin{equation*}
        \int_{ L } p(\f x,t)\, d L(\f x)
\end{equation*}
on a line $L$ \cite{BurEtAl05,PNHB07}. In practical experiments they
can be realized by a thin laser beam that is either part of
Fabry--Perot \cite{BurEtAl05} or Mach--Zehnder \cite{PNHB07}
interferometer. Details on the technical realization of line
detectors can be found in \cite{BurEtAl05,PNHB07}. We just note here
that optical sensors offer high bandwidth and high signal to noise
ratio.


\subsection*{Data acquisition}

Throughout this paper we assume that the initial data $f$ is
supported in the open set $\Vrec$. We call $V_{\rm rec}$  {\em
recording volume}, if $\Vrec$ is invariant with respect to rotations
around the
          $\f e_3 = (0,0,1)$ axis, i.e.
          \begin{equation*}
                \Vrec = \set{R \f x: \f x \in \Vrec}
          \end{equation*}
          holds for all rotations $R \in {\rm SO(3)}$ with fix point set  $\R \f
          e_3$.
Important examples for recording volumes are the half space $\R^2
\times (0,\infty)$, the cylinder $D \times \R$, where $D$ denotes
the unit disc in $\R^2$. We set $\Drec := \Vrec \cap ( \set{0}
\times \R^2)$.

The measured data is collected with an array of parallel line
detectors that are tangential to $\Vrec$, orthogonal to $\f e_3$ and
rotated around the axis $\f e_3$. The rotation $R_{\fs \si}$ around
the axis $\f e_3$ is parameterized by $\fs \si \in S^1$ and defined
by
\begin{equation*}
    \begin{aligned}
        R_{\fs \si} \f e_1  &= \f e_1(\fs \si) := (\si_1, \si_2, 0)\,,\\
        R_{\fs \si} \f e_2  &= \f e_2(\fs \si) := (-\si_2, \si_1, 0)\,,\\
        R_{\fs \si} \f e_3  &= \f e_3\;.
    \end{aligned}
\end{equation*}
Here $(\f e_1,\f e_2,\f e_3)$ denotes the standard basis of $\R^3$
and $(\f e_1(\fs \si),\f e_2(\fs \si),\f e_3)$ is a positively
oriented orthonormal basis. In practical applications the
measurements, for fixed  $\fs \si$, are performed  either with a
array or by moving a single line detector along a recording curve
$\Crec$.

Let
\begin{equation*}
        L(\fs \si, \f y) := \set{ s \f e_1(\fs \si) +
                                  y_1 \f e_2(\fs \si) +
                                  y_2 \f e_3  : s \in \R }
\end{equation*}
be a line, where $\f y = (y_1, y_2) \in \R^2$. The positions of the
line detectors are given by $L(\fs \si, \f y_{\rm rec})$ with $\f
y_{\rm rec} \in \Crec \tm \partial \Drec$ which is called {\em
recording curve}, see Figure \ref{fg:line:geo}.

Since $\Vrec$ is rotationally invariant with respect to $\f e_{3}$
and $\Crec \tm \partial \Drec$ all lines $L(\fs \si, \f y_{\rm rec})$
(locations of the line detectors) are outside the support of $f$.
The forward operator
\begin{equation*}
    \begin{aligned}
        \Po: C_0^\infty( \Vrec )  & \to  C^\infty \left( S^1 \times  \Crec  \times (0, \infty) \right) \\
                        f               & \mapsto  \Po(f) := \left( (\fs \si, \f y_{\rm rec}, t) \mapsto
                                                             \int_{L(\fs \si, \f y_{\rm rec})} p(\f x,t) \, d L(\f x) \right)
    \end{aligned}
\end{equation*}
collects the data measured by all line detectors $L(\fs \si,\f
y_{\rm rec})$ for given initial data $f$. Here $p(\f x,t)$ denotes
the unique solution of \eqref{eq:wave3d}, \eqref{eq:ini3d1},
\eqref{eq:ini3d2}. We note that \eqref{eq:wave3d},
\eqref{eq:ini3d1}, \eqref{eq:ini3d2} is a well posed problem and
therefore $\Po(f)$ can be calculated stable for arbitrary initial
data $f$. Thermoacoustic CT with line detectors deals with the
solution of the operator equation $\Po( f ) = g $ with given
(potentially  noisy) measured data $g$. Therefore we study
uniqueness, stability and explicit inversion of $\Po$. In this paper
we focus on the derivation and implementation of inversion formulas.

\begin{psfrags}
    \psfrag{Lrec}{$L(\fs \si, \f y_{\rm rec})$}
    \psfrag{Crec}{$\Crec$}
    \psfrag{E3}{$\f e_3$}
    \psfrag{E2}{$\f e_2$}
    \psfrag{S}{rotation ($\fs \si$)}
    \begin{figure}[htb]
        \begin{center}
            \includegraphics[width = 0.7\textwidth]{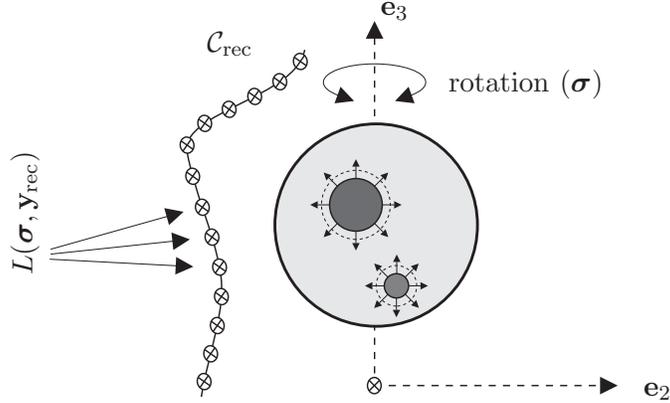}
        \end{center}
    \caption{{\bf Scanning geometry for thermoacoustic CT using line detectors.}
             The line detectors (orthogonal to the drawing plane) are uniformly
             distributed on the recording curve $\Crec$ and rotated around
             a single axis of rotation $\f e_3$, indicated by parameter $\fs \si$.}
    \label{fg:line:geo}
    \end{figure}
\end{psfrags}

\section{Decomposition of the three dimensional forward operator into a system
         of two dimensional operators}\label{sec:rot}

Let $V_{\rm rec}$ be a recording volume and $\Crec \tm \partial
\Drec$ a recording curve. In this Section we prove that the
corresponding forward operator $\Po$ decomposes into the product of
the two dimensional Radon transform $\Ro$ and the operator
\begin{equation*}
    \begin{aligned}
        \Qo : C_0^\infty( \Drec )   & \to  C^{\infty} \big( \Crec \times  (0, \infty) \big) \\
                                F         & \mapsto \Qo(F):= \big( (\f y_{\rm rec}, t) \mapsto  P( \f y_{\rm rec} , t) \big)
    \end{aligned}
\end{equation*}
that maps a planar function $ F \in C_0^\infty( \Drec ) $
onto the solution $P$ of the two dimensional wave equation
\begin{equation}\label{eq:wave2d}
    \left(  \frac{\partial^2}{\partial t^2}
    - \lap_{\f y} \right) P = 0 \,, \qquad
    (\f y, t)  \in \R^2  \times (0, \infty)
\end{equation}
with initial conditions
\begin{eqnarray}
  P (\f y , 0) =  F(\f y)\,,                            & \f y \in \R^2 \,, \label{eq:ini2d1} \\
  \frac{\partial P}{\partial t} (\f y , 0 ) = 0  \,,    & \f y \in \R^2 \label{eq:ini2d2}
\end{eqnarray}
restricted to the recording curve $\Crec$. Here $\lap_{\f y}$
denotes the Laplacian in $\R^2$. We use the following notation:
\begin{enumerate}
    \item Define
          \begin{equation*}
            \begin{aligned}
                \Ro \otimes \Io_3: C_{0}^{\infty}(\Vrec) & \to C_{0}^{\infty}(S^{1} \times \R^{2}) \\
                                                      f        & \mapsto (\Ro \otimes \Io_3)(f):=
                                                                          \Big( \big( \fs \tau, (r,z) \big) \mapsto
                                                                          \big( {\cal R} (f_{z}) \big)(\fs \tau, r) \Big)
            \end{aligned}
          \end{equation*}
          where $\Ro$ denotes the classical two dimensional Radon transform
          \cite{Nat86} defined by
          \begin{equation*}
            \begin{aligned}
                \Ro: C_{0}^{\infty}(\R^{2}) & \to C_{0}^{\infty}(S^{1} \times \R) \\
                                   \ph  & \mapsto \Ro(\ph):= \left( (\fs \tau, r) \mapsto \int_{\R} \ph(s \fs \tau + r \fs \tau^{\perp}) \, ds  \right)
            \end{aligned}
          \end{equation*}
          and $f_{z}$ the function $f$ for fixed third parameter $z$, i.e.
          \begin{equation*}
            f_{z}(x,y) := f(x,y,z)\,, \quad \textrm{with $z \in \R$ fixed\;.}
          \end{equation*}
          Here $\fs \tau^\bot$ is a unit vector orthogonal to
          $\fs \tau$.
          The notation indicates that the operator $\Ro \otimes \Io_3$ operates on
          a function $f$ defined on $\R^3$ by applying the Radon transform to the first two
          components for fixed third component.
    \item For a fixed $\fs \tau \in S^{1}$ and a function $G \in C_{0}^{\infty}(S^{1} \times \Drec)$ define
          \begin{equation*}
            G_{\fs \tau}(x_{1}, x_{2}) := G(\fs \tau, x_{1}, x_{2}) \;.
          \end{equation*}
          Furthermore, we define an operator $\Io \otimes \Qo$ by
          \begin{equation*}
            \begin{aligned}
                \Io \otimes \Qo: C_{0}^{\infty}(S^{1} \times \Drec) & \to C^{\infty} \big( S^{1} \times \Crec \times (0, \infty) \big) \\
                                                                G         & \mapsto (\Io \otimes \Qo)(G):= \big( (\fs \tau, \f y_{\rm rec},t)
                                                                            \mapsto \big( \Qo (G_{\fs \tau}) \big)(\f y_{\rm rec},t) \big) \;.
            \end{aligned}
          \end{equation*}
\end{enumerate}

The main feature of using line detectors is based on the following decomposition:

\begin{theorem}\label{thm:decomp}
    Let  $f \in C_0^\infty( \Vrec )$. Then
    \begin{equation}\label{eq:decomp}
        \Po (f)  = (\Io \otimes \Qo) \circ (\Ro \otimes \Io_3)(f)\;.
    \end{equation}
\end{theorem}

\begin{proof}
    Let $ f \in C^\infty_0 (\Vrec)$,  $p$ be the solution of \eqref{eq:wave3d},
    \eqref{eq:ini3d1}, \eqref{eq:ini3d2} and fix $\fs \si \in S^{1}$. We define
    \begin{equation*}
        P(\f y,t) := \int_{L(\fs \si,\f y)} p(\f x,t) \, dL(\f x)
    \end{equation*}
    and show that $P$ is the unique solution of the two dimensional wave
    equation \eqref{eq:wave2d}, \eqref{eq:ini2d1}, \eqref{eq:ini2d2} with initial data
    $F = (\Ro \otimes \Io_3)(f) (\fs \si, \cdot)$. Theorem \ref{thm:decomp} then
    follows from the fact that $\Po(f)(\fs \si, \f y_{\rm rec},t) = P(\f y_{\rm rec},t)$
    for $\f y_{\rm rec} \in \Crec$ and $t\geq 0$. From \eqref{eq:ini3d1},
    \eqref{eq:ini3d2} and the definitions of $(\Ro \otimes \Io_3)(f)$ and $P$ it follows
    immediately that $P$ satisfies the initial conditions \eqref{eq:ini2d1}, \eqref{eq:ini2d2}.
    It remains to prove that $P$ satisfies \eqref{eq:wave2d}.

    Since $(\f e_1(\fs \si), \f e_2(\fs \si), \f e_3)$ is an orthonormal
    basis of $\R^3$ for every $\f x \in \R^3$ there exist unique elements $s \in \R$ and
    $\f y \in \R^2$ such that $\f x = s \f e_1( \fs \si) + \f y =: s \f e_1(\fs \si)
    + y_1 \f e_2(\fs \si) + y_2 \f e_3$.
    Therefore, the Laplacian in $\R^3$ decomposes into the sum
    \begin{equation*}
        \lap_{\f x} = \frac{\partial^2 }{\partial s^2} + \lap_{\f y}\;.
    \end{equation*}
    Using \eqref{eq:wave3d} results in
    \begin{equation*}
        \begin{aligned}
            0 &= \int_{\R} \left(\frac{\partial^2 p}{\partial t^2}(s \f e_1(\fs \si) + \f y, t)
                 - \frac{\partial^2 p}{\partial s^2}(s \f e_1(\fs \si) + \f y, t)
                 - \lap_{\f y} p (s \f e_1( \fs \si) + \f y, t) \right) ds
                 \\
              &= \frac{\partial^2}{\partial t^2} \int_{\R} p ( s \f e_1(\fs \si) + \f y, t) \, ds
                 - \left[ \frac{\partial p}{\partial s}(s \f e_1(\fs \si) + \f y, t) \right]_{s = -\infty}^{s = \infty}
                 \\
              &  \phantom{=~} - \lap_{\f y} \int_{\R} p(s \f e_1(\fs \si) + \f y, t) \, ds
                 \\
              &= \frac{\partial^2 P}{\partial t^2}(\f y,t) - \lap_{\f y} P(\f y, t) \;.
        \end{aligned}
    \end{equation*}
    The last equality holds since $p(\cdot, t)$ is compactly supported \cite{Joh82}
    for all $t \geq 0$. Hence $P$ satisfies \eqref{eq:wave2d}.
\end{proof}

The decomposition of $\Po$ into a set of lower dimensional operators
can be used to reduces the complexity of the derived reconstruction
algorithms.

\medskip
The operator $\Qo$ is closely related to the circular Radon
transform with center set $\Crec$ \cite{FinPatRak04,Pal04}, which is
defined by
\begin{equation*}
    \begin{aligned}
        \Mo: C_0^\infty(\Drec) & \to C^\infty \big( \Crec \times (0, \infty) \big) \\
                            F        & \mapsto \Mo (F) :=  \left( (\f y_{\rm rec}, r) \mapsto \frac{1}{2 \pi}
                                                           \int_{S^{1}} F( \f y_{\rm rec} + r \fs \om) \,d \Om(\fs \om) \right)
    \end{aligned}
\end{equation*}
and maps a function $F \in C_0^\infty ( \Drec )$ onto its integrals
over circles with centers on the recording curve $\Crec$.

Using D'Alembert's formula the solution of \eqref{eq:wave2d},
\eqref{eq:ini2d1}, \eqref{eq:ini2d2} turns out to be \cite{Joh82}
\begin{equation}\label{eq:dalambert2D}
    \big(\Qo(F) \big) (\f y_{\rm rec}, t) = \frac{\partial}{\partial t} \int_{0}^{t}
                                            \frac{r \Mo (F)(\f y_{\rm rec},r)}{\sqrt{t^2 - r^2}} \, dr
\end{equation}
which means that $\Qo(F)$ can be calculated if the circular Radon
transform $\Mo(F)$ of $F$ is known. The following Lemma states that
also the reversal is true:
\begin{lemma} \label{lemma:qm}
    Let $F \in C_0^\infty( \Drec )$. Then
    \begin{equation}\label{eq:findM}
        \Mo (F) (\f y_{\rm rec}, r) = \frac{2}{\pi}  \int_{0}^{r} \frac{\Qo (F)(\f y_{\rm rec},t)}{\sqrt{r^2 - t^2}} \; dt
    \end{equation}
    for all $\f y_{\rm rec} \in \Crec$ and $t \geq 0$.
\end{lemma}

A proof of a generalization of Lemma \ref{lemma:qm} is given in
\cite[Theorem 2.1.2]{Tri01}, but for the sake of completeness we
prove Lemma \ref{lemma:qm} since it is an elementary proof.
\begin{proof}
    Let $F \in C_0^\infty(\Drec)$, fix $\f y_{\rm rec} \in \Crec$ and define
    $\Phi(t) := \Mo(F)(\f y_{\rm rec}, t)$. Furthermore, we define for a function
    $\Psi \in C^\infty \big( (0, \infty) \big)$
    \begin{equation*}
        \begin{aligned}
            \big( \Jo(\Psi) \big)(u)            &:= \int_{0}^{u} \frac{\Psi(t)}{\sqrt{u^2 - t^2}} \; dt \,,\\
            \big( \Jo_{\rm mul}(\Psi) \big)(u)  &:= \int_{0}^{u} t \frac{\Psi(t)}{\sqrt{u^2 - t^2}} \; dt \,, \\
            \big( \Ko(\Psi) \big)(u)            &:= u \big( \Jo(\Psi)\big)(u) \;.
        \end{aligned}
    \end{equation*}

    Due to equation \req{dalambert2D} it is sufficient to prove that
    $(2/ \pi) \Jo \big( \big( \Jo_{\rm mul}(\Phi) \big)' \big) = \Phi$.
    As an auxiliary result we prove that for all $\Psi \in C^\infty((0, \infty))$
    \begin{equation}\label{equ:idR}
        I(r) := \big( \Jo \big( \Ko (\Psi) \big) \big)(r) = \frac{\pi}{2} \int_0^r \Psi(u) du \;.
    \end{equation}
    For the proof of (\ref{equ:idR}) let $t > 0$ and consider the integral
    \begin{equation*}
        I(r) =  \int_0^r  \frac{t\Jo(\Psi)(t)}{\sqrt{r^2-t^2}} \; dt
             =  \int_0^r \left( \int_0^{t} \Psi(u) \frac{du}{\sqrt{t^2-u^2}}\right) \frac{t}{\sqrt{r^2-t^2}} \; dt \;.
    \end{equation*}
    From {\em Fubini's theorem} we conclude that
    \begin{equation*}
        I(r) = \int_0^r \left( \int_u^r \frac{t}{\sqrt{t^2-u^2} \sqrt{r^2-t^2}} \; dt \right) \Psi(u) du
    \end{equation*}
    We substitute $t = \sqrt{u^2 + v^2}$ and $dt =  (v/t) dv$ in the inner integral. Therefore
    \begin{equation*}
        I(r) = \int_0^r \left( \int_0^{\sqrt{r^2 - u^2}}
               \frac{dv}{\sqrt{r^2 - u^2 - v^2}} \right) \Psi(u) du \;.
    \end{equation*}
    Finally, we use the substitution $v = \sqrt{r^2 - u^2} \sin(\al)$ and $dv = \sqrt{r^2-u^2} \cos(\al) d\al$.
    We obtain
    \begin{equation*}
    I(r) = \int_0^r \left( \int_0^{\pi/2} d\al \right) \Psi(u) du
         = \frac{\pi}{2} \int_0^r \Psi(u) du \;.
    \end{equation*}
    Hence we have proved (\ref{equ:idR}).

    Next we prove that
    \begin{equation}\label{equ:bo}
        \big( \Jo_{\rm mul}(\Phi) \big)'  = \Ko(\Phi') \;.
    \end{equation}
    Using the relation
    \begin{equation*}
        \frac{\partial}{\partial r} \sqrt{t^2 - r^2} = - \frac{r}{\sqrt{t^2 - r^2}}
    \end{equation*}
    and integration by parts we obtain
    \begin{equation*}
        \begin{aligned}
            \int_{0}^{t} \frac{r \Phi(r)}{\sqrt{t^2 - r^2}} dr
            &=
            \left[ \Phi(r) \sqrt{t^2 - r^2} \right]_{r=0}^{r=t} + \int_{0}^{t} \Phi'(r) {\sqrt{t^2 - r^2}} dr
            \\
            &=
            \int_{0}^{t} \Phi'(r) {\sqrt{t^2 - r^2}} dr \;.
        \end{aligned}
    \end{equation*}
    By differentiating the above equation  with respect to $t$ we obtain
    \begin{equation*}
        \begin{aligned}
            \left( \int_{0}^{t}  \frac{r \Phi(r)}{\sqrt{t^2 - r^2}} dr \right)'(t)
            &= \left[ \Phi'(r) \sqrt{t^2 - r^2} \right]_{r=t} + t \int_{0}^{t} \frac{\Phi'(r)}{\sqrt{t^2 - r^2}}\, dr \;.\\
            &= \big( \Ko (\Phi') \big)(t)
        \end{aligned}
    \end{equation*}
    Hence we have proved \eqref{equ:bo}.
    Applying $\Jo$ to \eqref{equ:bo} we find
    \begin{equation*}
        \Jo \big( \big( \Jo_{\rm mul}(\Phi) \big)' \big)  =  \Jo \big( \Ko(\Phi') \big) = \frac{\pi}{2} \Phi \;.
    \end{equation*}
For the last equality we applied the auxiliary result
(\ref{equ:idR}) to $\Psi = \Phi'$.
\end{proof}

In order to prove that $\Po$ is injective we use Lemma
\ref{lemma:qm} and a uniqueness result for the circular  Radon transform.

\begin{theorem}[Uniqueness for arbitrary recording curve]\label{cor:unique}
Let $\Vrec$, $\Drec$ and $\Crec$ be as above, and assume
additionally, that $\Drec$ is convex and that $\Crec \subset \Drec$
is a relatively open $C^1$ curve. If $f \in C_0^\infty(\Vrec)$ and $
\Po (f) = 0$, then $f= 0$.
\end{theorem}

\begin{proof}
Due to  \req{decomp} it is sufficient to show that $\Ro$ and $\Qo$
are injective. It is well known that the Radon transform
$\Ro$ is injective, see for example \cite[Theorem 2]{Nat86}.
Therefore, it remains to show that $\Qo$ is injective.

Let $F \in C_0^\infty(\Drec)$ and $ \Qo (F)  = 0$. From \req{findM}
it follows that $\Mo(F) = 0$. It is known that the circular Radon
transform with center set $\Drec$ is injective, see \cite[Theorem
19]{art05:AmbKuc} (compare also with \cite{AgrQui96,FinPatRak04}).
We conclude that  $ F =0$ which means that $\Qo$ is injective.
\end{proof}

Equation \req{decomp} implies that thermoacoustic CT with line
detectors deals with the inversion of the operator $\Qo$ and the
inversion of the two dimensional classical Radon transform $\Ro$.
Fast and stable algorithms for inverting the two dimensional
classical Radon transform $\Ro$ are well investigated. The key task
for deriving inversion algorithms for $\Po$ is the inversion of the
operator $\Qo$. Therefore, in the next Section we  present  and
analyze an analytic formula for inverting $ \Qo $.

\section{Linear recording curve: Inversion of the two dimensional wave equation
from data on a line}
\label{sec:linear}

\begin{figure}[h!]
    \begin{psfrags}
        \psfrag{S}{$\fs \si$}
        \psfrag{E3}{$\f e_{3}$}
        \psfrag{Crec}{$\Crec$}
        \begin{center}
            \includegraphics[height=0.45\textwidth]{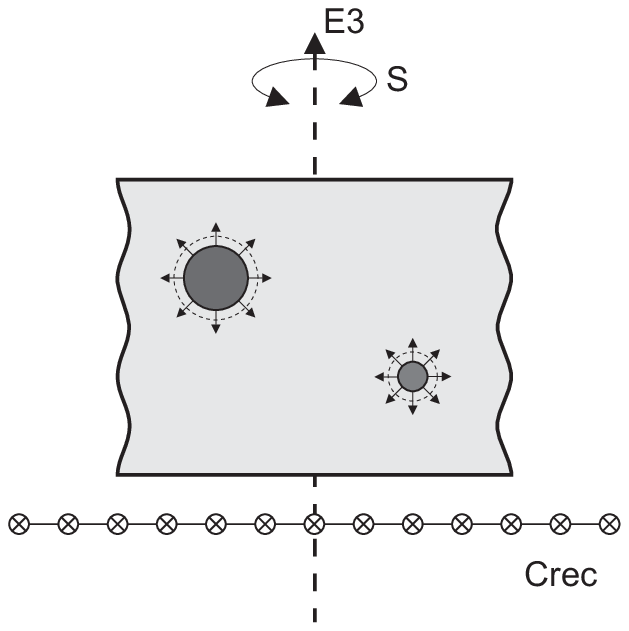}\hspace{7em}
            \includegraphics[height=0.4\textwidth]{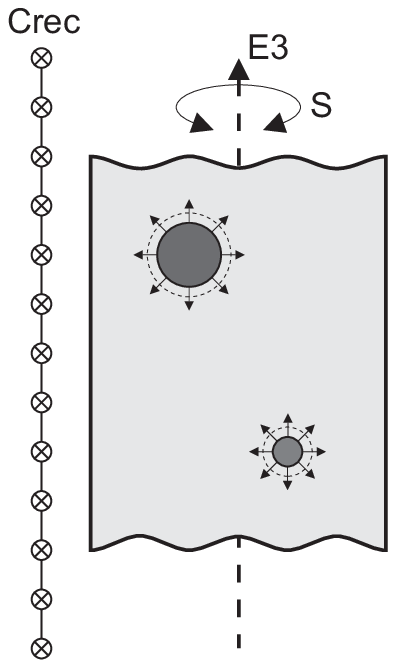}\\[1em]
            \includegraphics[height=0.45\textwidth]{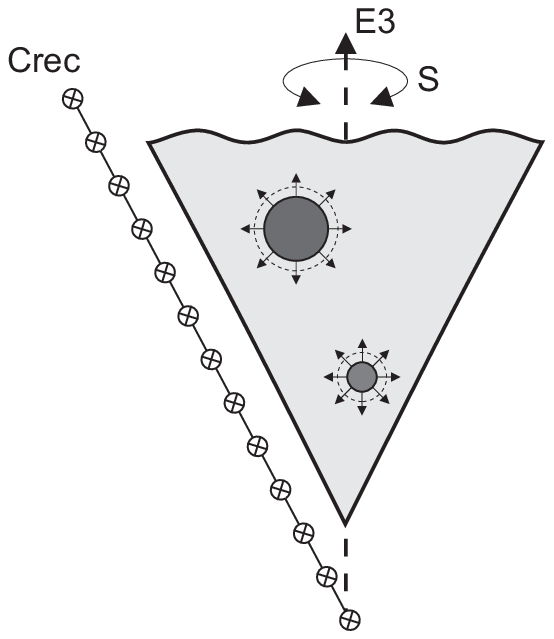}\hspace{7em}
            \includegraphics[height=0.45\textwidth]{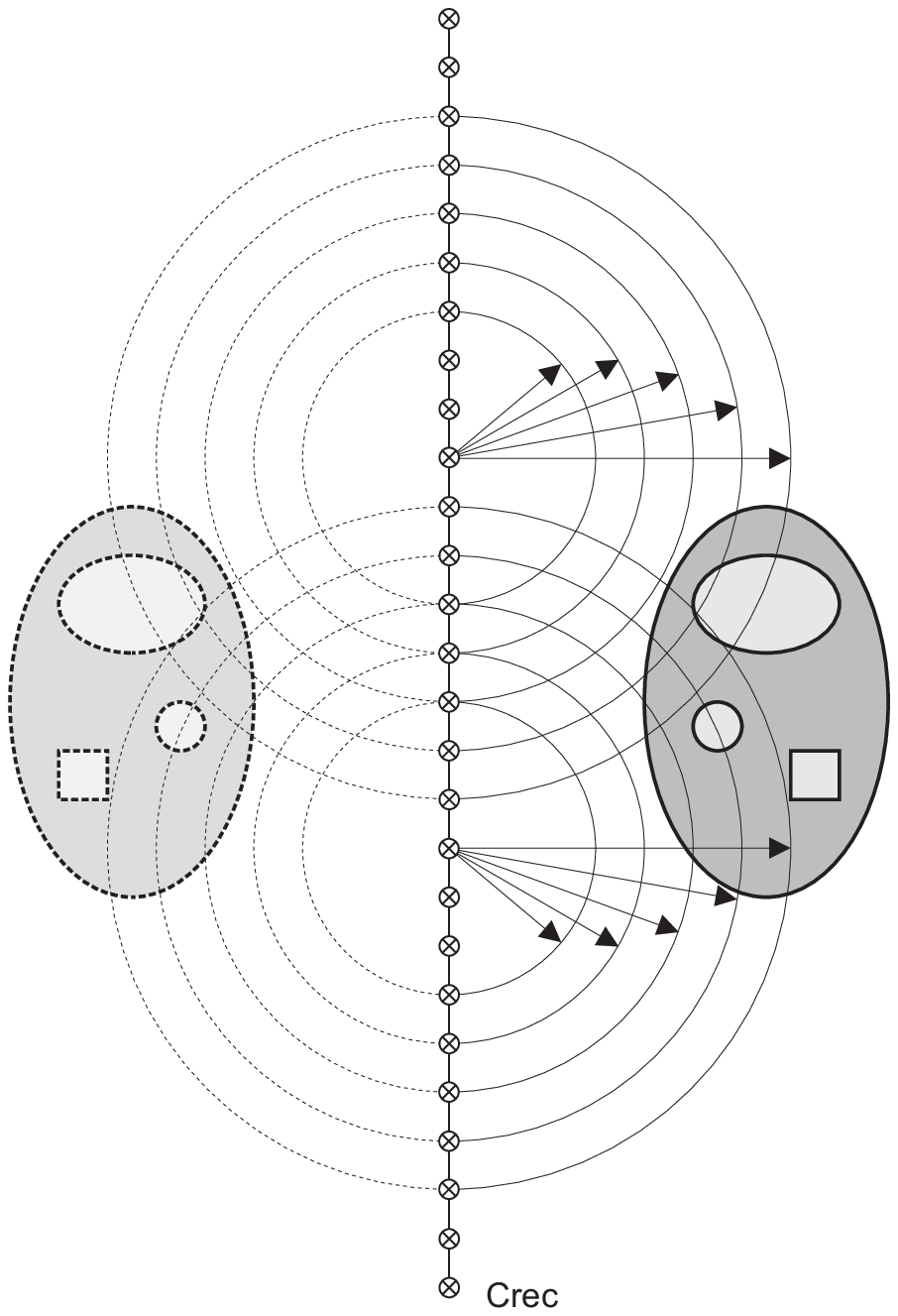}\\
        \end{center}
        \caption{
                The array of line detectors measures the solution of the two dimensional
                wave equation \req{wave2d}-\req{ini2d2}
                restricted to $\Crec$ and is rotated in a plane   (top left),
                tangential to a cylinder (top right) or  tangential to a cone (bottom left).
                Bottom right: The restriction of the solution to
                $\Crec$ is uniquely determined by the means over circles with centers on
                $\Crec$ (see \req{dalambert2D}) and remains unchanged after symmetrizing
                the initial data $F$ around $\Crec$.
                }
    \label{fg:linear12}
    \end{psfrags}
\end{figure}

Throughout the remaining parts of this article we consider the case
where the recording curve  $\Crec \tm \R^2$ is a line. Depending on
the angle  between $\Crec$ and $\f e_2$, the recording volume
$V_{\rm rec}$ is either a closed cylinder, a closed cone or a half
space, each considered as a subset of $\R^{3}$ (see Figure
\ref{fg:linear12}). In all cases $\Qo$  maps the initial data $F$
onto the solution of of the two dimensional wave equation
\req{wave2d}-\req{ini2d2} restricted to $\Crec$.

From symmetry properties of the wave equation (bottom right image in
Figure  \ref{fg:linear12}) it follows that in order to invert $\Qo$
it is sufficient to solve the following problem: Recover $H \in
C_0^\infty(\RtI)$ from data $G := P|_{\R \times \set{0} \times (0,
\infty)}$, where $P$ is  the solution of
\begin{equation}\label{eq:wave-sym}
    \left(  \frac{\partial^2}{\partial t^2}
    - \frac{\partial^2}{\partial u^2} - \frac{\partial^2}{\partial v^2  } \right) P = 0 \,,\quad
    (u, v, t)  \in \R^2  \times (0, \infty)
\end{equation}
with symmetric  initial conditions
\begin{eqnarray}
  P (  u,  v , 0) =  H(u , \abs{v}  )\,,    & (u, v) \in \R^2 \,, \label{eq:ini-sym1} \\
  \frac{\partial P}{\partial t} (u, v, 0 ) = 0  \,, & (u, v ) \in \R^2 \label{eq:ini-sym2}.
\end{eqnarray}
$u$ denotes the coordinate in direction $\Crec$ and $v$ the signed
distance to the orthogonal direction.

In the following we have to apply the Fourier--cosine transform to
$G$. Since $G \not\in L^1(\RtI)$, see Corollary \ref{cor:notL1} below, we have to sidestep to the
space of symmetric tempered distributions. See Appendix
\ref{sec:app} for the definitions of the spaces $\Ssym(\RtR)$,
$\Czzsym(\RtR)$, $\Spsym(\RtR)$, the Fourier--cosine
transforms $\cost$, $\costp$ and the
natural injections $\injS$, $\injL$, $\injC $ that map functions
onto tempered distributions. We note that, since $G_{\rm sym} \in
\Czz(\RtR)$, $\injC[G_{\rm sym}] \in \Spsym(\RtR)$
and the Fourier--cosine transform $\costp [ \injC [G_{\rm sym}] ]$
is well defined, see Appendix \ref{sec:app}
for details.

\begin{theorem}[Fourier inversion formula] \label{thm:fif}
Let $H \in C_0^\infty(\RtI)$,  $P$ be the solution of
\req{wave-sym}-\req{ini-sym2}, $G := P|_{\R \times \set{0} \times
(0, \infty)}$ and $\Om := \set{(\om,\la) \in \R^2 : \om >
\abs{\la}}$. Then there exists $\bar{G} \in  L^1(\RtI) \cap
C^{\infty}(\Om) $ with $\costp [ \injC[  G_{\rm sym}  ] ]  = \injL[
\bar{G}_{\rm sym} ]$, and the Fourier inversion formula holds
\begin{equation}\label{eq:fif}
    \cost [ H ] ( \la, \ka  ) =
                \bar{G} \left(  \la, \sqrt{ \ka^2 + \la^2 } \right)
                \cdot \displaystyle \frac{\ka}{\sqrt{ \ka^2 + \la^2  }}
\end{equation}
point--wise for all $(\la, \ka)  \in \R \times (0,\infty)$.
\end{theorem}

\begin{proof}
An elementary calculation shows that for all $(\la, \ka) \in \RtI$
the functions
\begin{equation*}
    \cos(\ka v) e^{i \la u} \cos \left( t \sqrt{\ka^2 + \la^2} \right)
\end{equation*}
are solutions of (\ref{eq:wave-sym}) and (\ref{eq:ini-sym2}).
Since $H \in C_0^\infty( \RtI) \tm  \Ssym(\RtR)$,
$\bar{H}  := \cost[H]$ is well defined.
From the linearity of the wave equation and the
inversion formula for the Fourier--cosine transform
on $\Ssym(\RtR)$, see Proposition \ref{proposition:iso-s} in Appendix \ref{sec:app},
it follows that the unique solution of (\ref{eq:wave-sym})-(\ref{eq:ini-sym2}) is
given by
\begin{equation*}
    P (u, v, t) =   \frac{1}{\pi}
                    \int_{\R}
                    \left(
                           \int_0^{\infty} \bar{H}( \la, \ka )
                           \cos \left(t \sqrt{\ka^2 + \la^2} \right)
                           \cos(\ka v) \, d \ka
                    \right)
                    e^{i \la u }  d \la \;.
\end{equation*}
Substituting $ \ka^2 = \om^2 - \la^2 $ in the inner integral
and taking $v = 0$ leads to
\begin{equation}\label{eq:Qft}
    G(u,  t) =  \frac{1}{\pi}
                \int_{\R}
                \left(
                        \int_{\abs{\la}}^{\infty} \bar{H}
                        \left( \sqrt{\om^2 - \la^2},  \la \right)
                        \cos(\om t) \frac{ \om  d\om }{ \sqrt{ \om^2 - \la^2 }}
                \right)
                e^{i \la u }  d \la \;
\end{equation}
Define a function $\bar{G}$ by
\begin{equation}\label{eq:hatq}
    \bar{G} ( \la , \om ) := \begin{cases}
                                        \bar{H} ( \la, \sqrt{\om^2 - \la^2}  )
                                        \om / \sqrt{ \om^2 - \la^2 }  & \textrm{if } \om > \abs{\la} \\
                                        0 & \textrm{otherwise} \;.
                                    \end{cases}
\end{equation}
$\bar{G}$ satisfies
\[
    \norm{\bar{G}}_{L^1} =
    \int_R \int_0^\infty \abs{\bar{G}(\la, \om)} \, d\om \, d\la
    =
    \int_R \left(
    \int_{k}^{\infty} \abs{\bar{H} \left( \la, \sqrt{\om^2 - \la^2}   \right)}
    \frac{ \om d\om }{ \sqrt{ \om^2 - \la^2 }}
    \right) d \la  \;.
\]
Substituting $ \om^2 = \ka^2 + \la^2$  in the inner
integral we conclude that $ \norm{\bar{G}}_{L^1} =
\norm{\bar{H}}_{L^1}$ and hence $\bar{G}  \in L^1(\RtI)$.
Equation (\ref{eq:Qft}) and $\bar{G}  \in L^1(\RtI)$
imply that $\cost [\bar{G} ]^{-} = G $ and since the
definitions of the Fourier--cosine transforms on $S'_{\rm sym}( \RtR)$
and $L^1(\RtI)$ are compatible, see Proposition \ref{proposition:C-comp}
\begin{equation*}
    \icostp [ \injL[ \bar{G}_{\rm sym}]] = \injC[\cost[\bar{G}]_{\rm sym}^{-} ]
                                         = \injC[G_{\rm sym} ] \;.
\end{equation*}
Applying the inversion formula for the Fourier--cosine transform on
$S'_{\rm sym}(\RtR)$, see Proposition \ref{proposition:iso-s-str},
shows $\injL[ \bar{G}_{\rm sym} ] = \costp [ \injC[ G_{\rm sym} ] ] $.

Finally, solving \req{hatq} for $ \bar{H} $ shows that \req{fif} holds point--wise
for all $\ka=\sqrt{\om^2-\la^2} > 0$.
\end{proof}

\begin{corollary}\label{cor:notL1}
Let $H$, $G$ and $\bar{G}$ be as in Theorem \ref{thm:fif}. Moreover, assume that
$H \neq 0$ is  non--negative. Then $G$ is not absolutely integrable.
\end{corollary}

\begin{proof}
Assuming $G\in L^1(\RtI)$ it follows from Proposition \ref{proposition:iso-s} in Appendix \ref{sec:app}
that $\bar{G} = \cost[G] \in \Czz(\RtI)$.
Since $H$  is non--vanishing it follows that $ \cost[H] ( 0, 0 ) > 0$
and since $\cost[H]$ is continuous there exists $\eps > 0$, such that
$\cost[H] ( \eps, \eps ) > 0$.  Inserting the
sequence $(\eps+1/n ,  \eps + 2/n )$ in \req{hatq}  shows
that $\bar{G}$ is unbounded, contradicting  $ \bar{G} \in \Czz(\RtI)$.
\end{proof}

Corollary \ref{cor:notL1} implies that the integral
\begin{equation}\label{eq:CG-not}
   \frac{1}{\pi}
   \int_\R \int_0^\infty   G(u, t)\cos(\om t) e^{ -i \la u} dt\,  du
\end{equation}
is not absolutely convergent and therefore cannot be used directly
to find an analytic representation of $\bar{G}$. However, in the
following we show that \req{CG-not} exists in some generalized sense
and can be used to find an analytic representation of $\bar{G}$.

In the following
\begin{equation*}
    \begin{aligned}
        \ft_{1}[\Phi](\la)  &:=   \left( \frac{1}{2 \pi} \right)^{1/2}
                                  \int_\R
                                  \Phi(u) e^{- i \la u }    du \,,\quad \la \in \R \,,
        \\
        \cost_2 [\Phi](\om) &:=   \left( \frac{2}{\pi} \right)^{1/2}
                                  \int_0^\infty
                                  \Phi(t) \cos(\om t)   dt \,,\quad \om  > 0
    \end{aligned}
\end{equation*}
denote the one dimensional {\em Fourier transform on $L^1(\R)$} and
the {\em Fourier-cosine transform on $L^1((0,\infty))$},
respectively. When applied to a function defined on $\R \times
(0,\infty)$ $\ft_{1}$ acts on the first and $\cost_2$ on the second
component. Finally, we recall the following special case of Jordan's
theorem: If $\Phi  \in L^1((0, \infty) )$  is continuous
differentiable on an interval including the point $\om >0$, then
\begin{equation}\label{eq:jordan2}
   \Phi( \om )
    =
    \lim_{ K \to  \infty}
    \left( \frac{2}{\pi} \right)^{1/2} \int_{0}^{K}
     \cos( \om t) \cost_2 [\Phi ] ( t ) \, d t\;.
\end{equation}

This enables us to prove our final result, namely an explicit
expression for the calculation of $\bar{G}$.

\begin{theorem}[Analytic formula for $\bar{G}$] \label{thm:fifComp}
Let $H$, $G$, $\bar{G}$ and $\Om$ as in Theorem
\ref{thm:fif}. Then, for $(\lambda, \omega) \in  \Omega$
\begin{equation}\label{eq:fctData}
    \bar{G}  (\lambda, \omega ) =
    \lim_{K \to \infty}
    \frac{1}{\pi}  \int_0^K
    \left(
    \int_\R G  ( u, t )
     e^{ -i \lambda u}  du \right)
     \cos( \omega t )  \, dt \;.
\end{equation}
\end{theorem}

\begin{proof}
From Theorem \ref{thm:fif} it follows that $\bar G \in L^1(\R \times (0,\infty))$
and $\cost [\bar G] = G^-$. Fubini's theorem implies that
\begin{equation} \label{eq:help1}
    \ft_1 \circ \cost_2 [\bar G] = \cost [\bar G]  = G^-.
\end{equation}
Since the initial data $H$ of the two dimensional wave equation
\req{wave-sym}-\req{ini-sym2} is compactly supported,
$G(\cdot, t ) \in C_0^\infty(\R)$ for fixed $t \geq 0$ and therefore $\ft_1 [G] (\cdot,
t ) \in S(\R)$. Hence, applying $\ft_1$ to \req{help1} yields
\begin{equation} \label{eq:help2}
    \ft_{1}[G]^- = \cost_2[\bar{G}]^{-}
\end{equation}
Next we apply \req{jordan2} to $\bar{G}(\la, \cdot )$, for fixed
$\la$,
\begin{equation} \label{eq:JFixLa}
   \bar{G} ( \la,  \om )
    =
    \lim_{ K \to  \infty}
    \left( \frac{2}{\pi} \right)^{1/2}
    \int_{0}^{K}
     \cost_2 [ \bar{G}] (\la,  t )  \cos( \om t)  \, d t \;.
\end{equation}

In order to finish the proof we apply \req{help2} to the integrand
of \eqref{eq:JFixLa} and obtain
\begin{eqnarray*}
   \bar{G} ( \la,  \om )
   & = &
   \lim_{ K \to  \infty}
    \left( \frac{2}{\pi} \right)^{1/2}
    \int_{0}^{K}
     \ft_1 [ G ] ( \la,  t )  \cos( \om t)  \, d t
   \\
   & = &
    \lim_{ K \to  \infty}
    \frac{1}{\pi}  \int_{0}^{K}
     \left(  \int_{\R}  G(u,t)
     e^{ -i \la u  } du \right)  \cos( \om t )  \, d t \;.
\end{eqnarray*}
\end{proof}

Theorem \ref{thm:fifComp} shows that the calculation of the Fourier--cosine
transform in the space $\Spsym(\RtR)$ of tempered distributions can be avoided.
Moreover, in practical applications data  $G(u, t)$ can be collected on a
finite space-time domain only. The fact  that $G(\cdot, t)$ is compactly supported for
$t$ fixed and Theorem \ref{thm:fifComp} justify that, in order to approximately
calculate $\bar{G}$, the domain of integration in \req{fctData} can be replaced
by this finite domain. However, data truncation will introduce error in the
reconstruction.

\begin{figure}[h!]
    \begin{psfrags}
        \begin{center}
            \includegraphics[height=0.4\textwidth,width=0.48\textwidth]{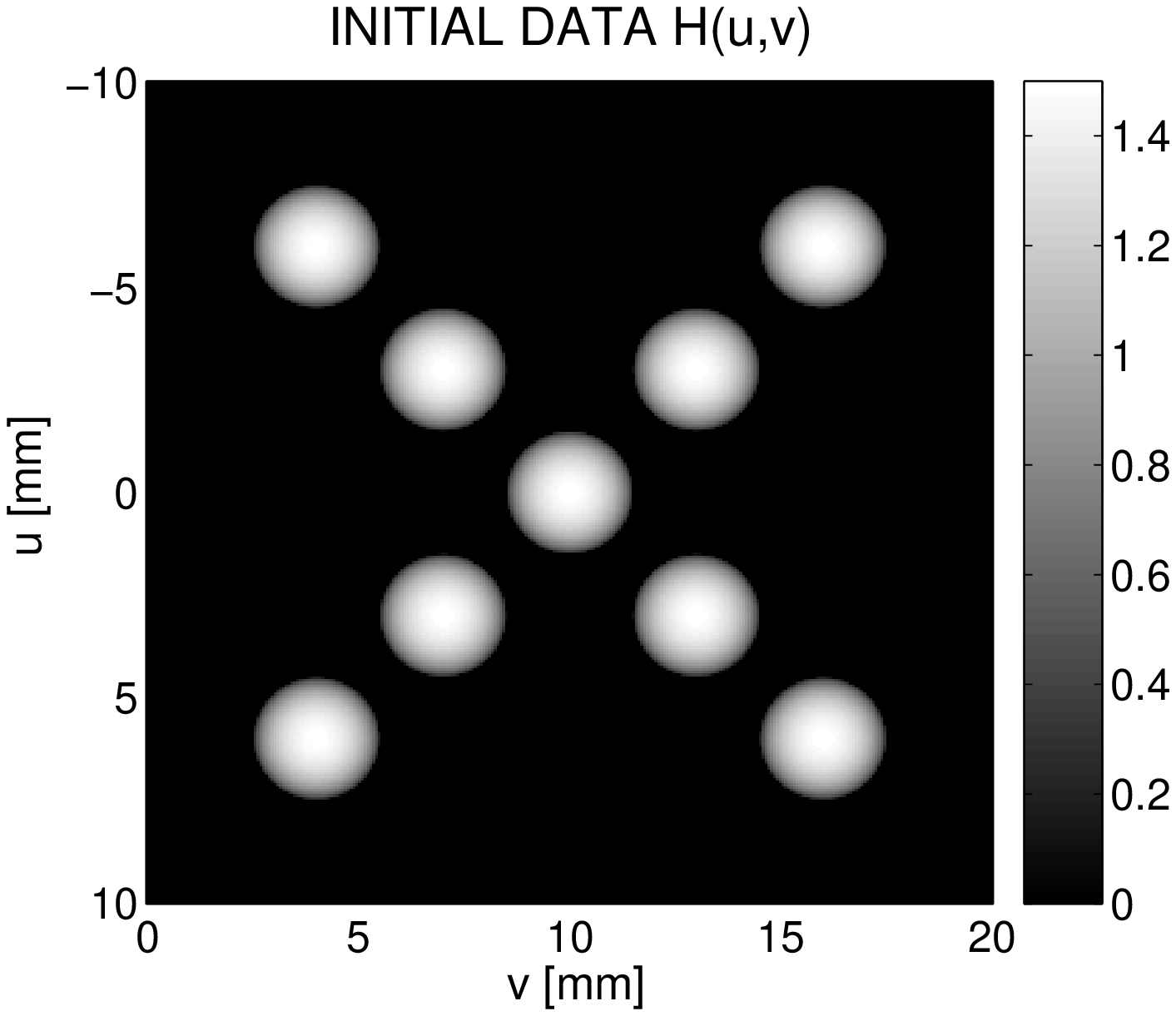}\hfill
            \includegraphics[height=0.4\textwidth,width=0.48\textwidth]{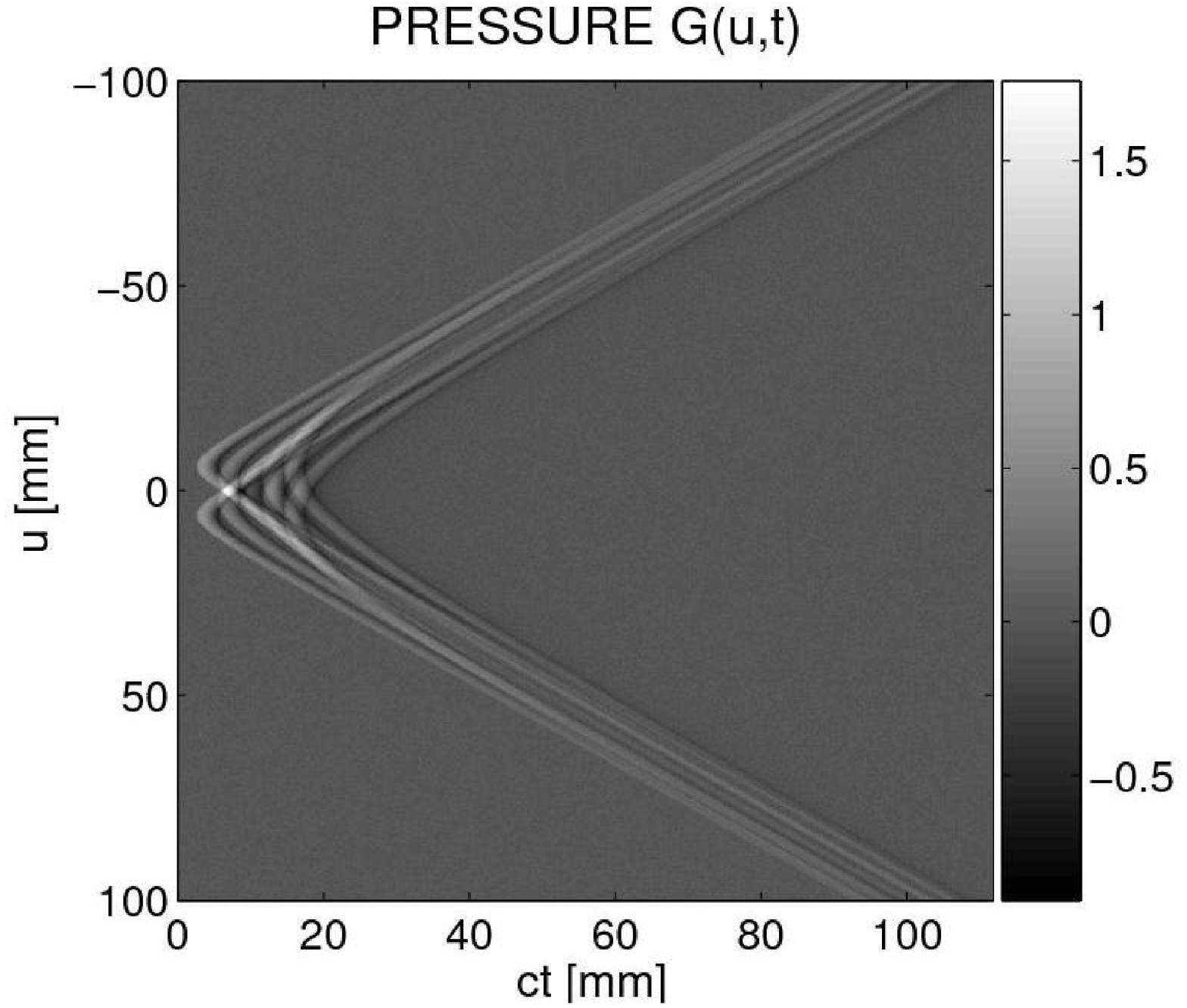}\\[1em]
            \includegraphics[height=0.4\textwidth,width=0.48\textwidth]{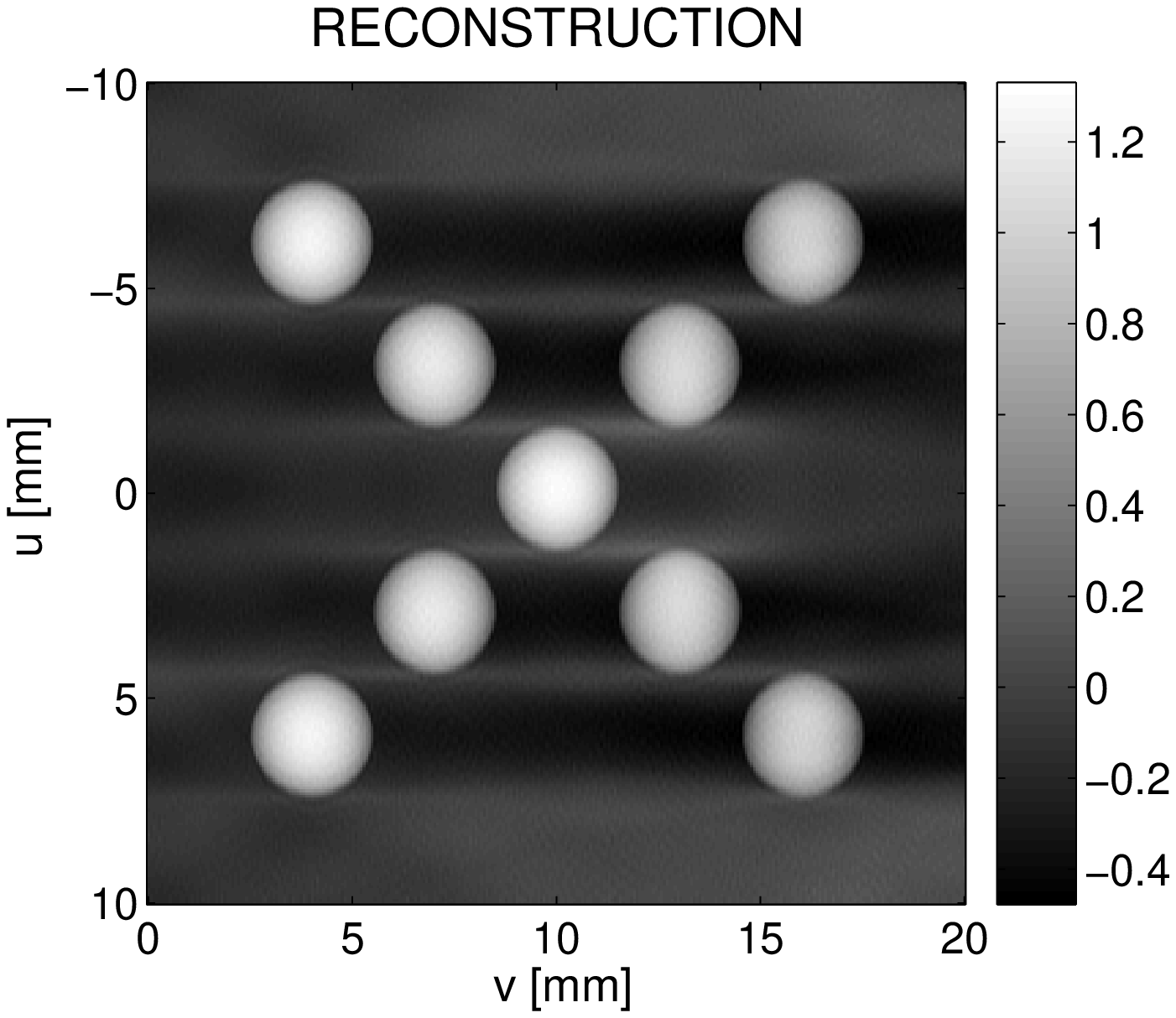}\hfill
            \includegraphics[height=0.4\textwidth,width=0.48\textwidth]{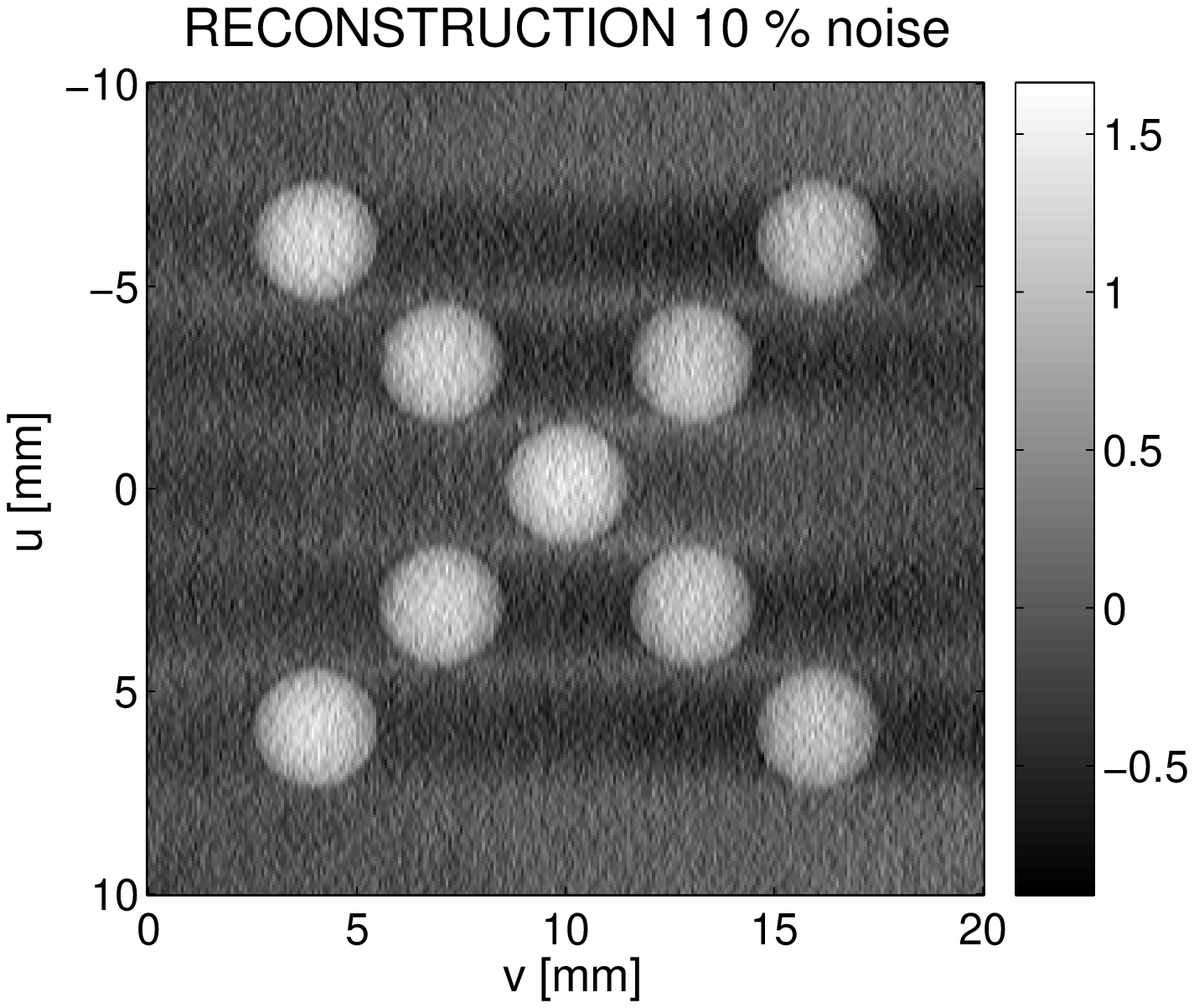}\\
        \end{center}
        \caption{
                {\bf Reconstruction with Fourier Formula.} Top left. Initial data $H$.
                Top right. Simulated data with  $ \abs{u} \leq 100$ mm.
                Bottom left. Reconstruction from exact data.
                Bottom right. Reconstruction from noisy data.}
    \label{fg:exfull}
    \end{psfrags}
\end{figure}

\begin{remark}[Numerical realization] \label{rem:num}
It is straight forward to derive a fast {\em Fourier reconstruction algorithm}
based on Theorems \ref{thm:fif}, \ref{thm:fifComp}:
First use the  {\tt FFT} algorithm to approximate $\bar{G}$,
see \req{fctData}. Then use bilinear interpolation and \req{fif} to approximate
$\cost [F]$ on a cartesian grid. Finally, again by using the {\tt FFT} algorithm,
calculate an approximation of $F$.
Due to the {\tt FFT} algorithm, using data $G$ to calculate $F$ at ${\tt N^2}$ grid points,
the numerical effort is $O( {\tt N^2} \log {\tt N})$ only.
\end{remark}

Figure \ref{fg:exfull} shows results of a numerical example with initial data
consisting of nine spheres supported in $ [-10,10] \times [0,20] $.
Data was collected with an array covering $[-100, 100]$.  The  solution of
(\ref{eq:wave-sym})-(\ref{eq:ini-sym2}) was calculated analytically and uniformly distributed
random noise with values in the interval $[-e,e]$, where $e$ equals
to $5 \%$ of the largest data value, was added. The bottom images are obtained via the
Fourier reconstruction formula, as presented in the Remark \ref{rem:num}.

\begin{figure}[h!]
    \begin{psfrags}
        \psfrag{u0}{$(u, 0)$}
        \psfrag{t1}{$\f t_1$}
        \psfrag{t2}{$\f t_2$}
        \psfrag{P1}{$P_1$}
        \psfrag{P2}{$P_2$}
        \begin{center}
            \includegraphics[height=0.4\textwidth,width=0.48\textwidth]{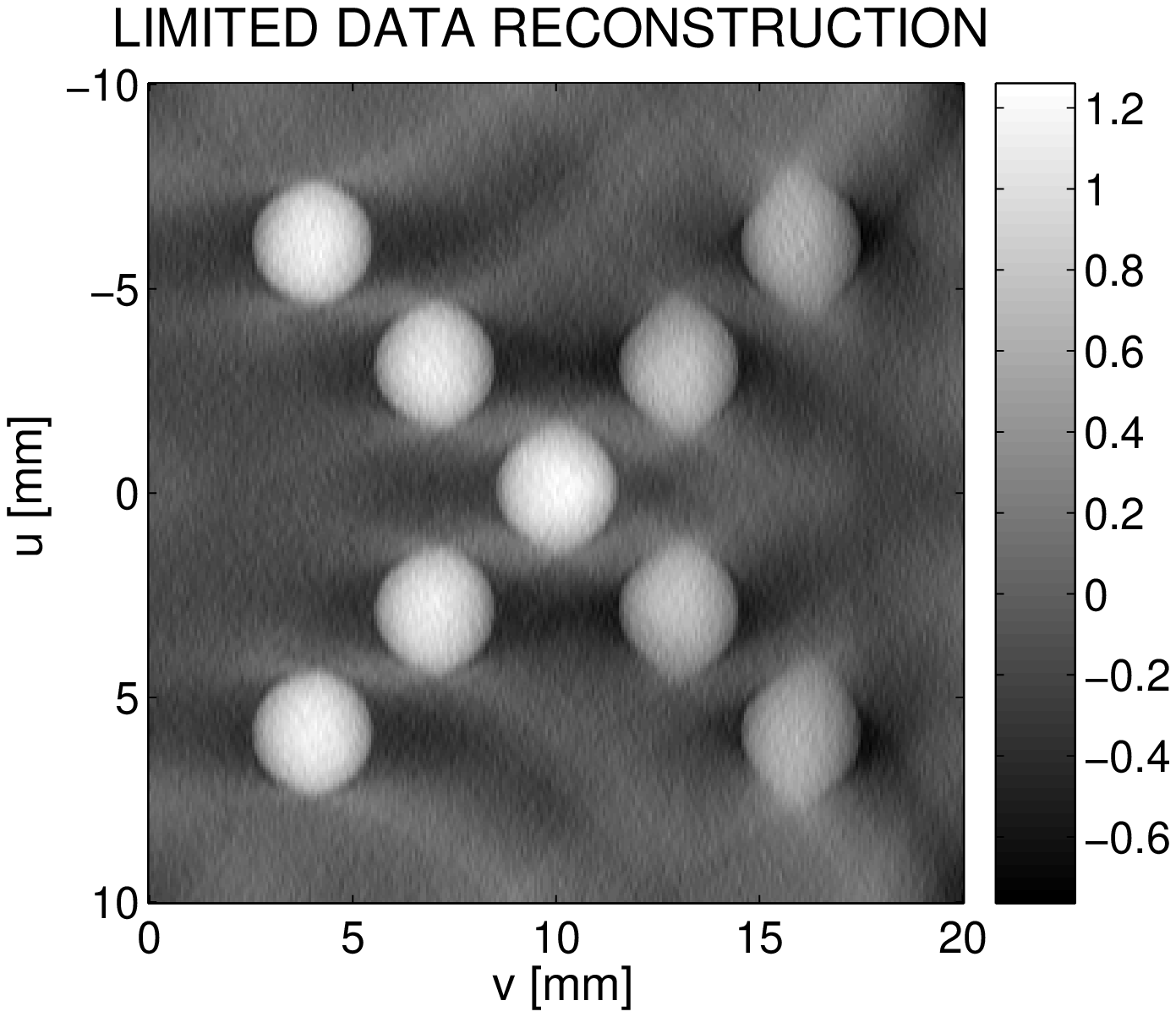}\hfill
            \includegraphics[width=0.48\textwidth]{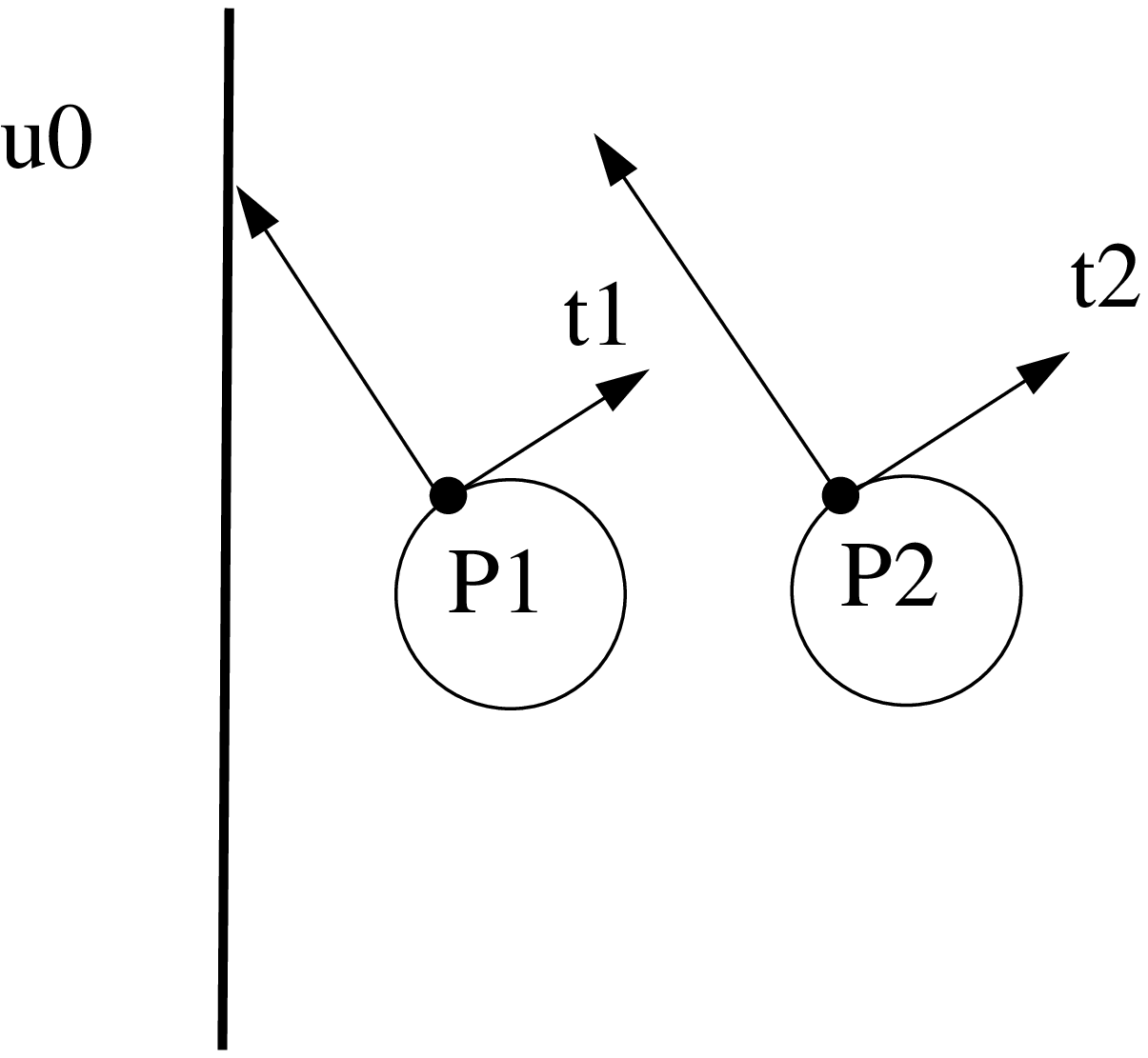} \\
        \end{center}
        \caption{ {\bf Finite aperture effect.}
            Left. Reconstruction from data $P(u, t)$ with $ \abs{u} \leq 20$ mm.
            Nearly horizontal boundaries are not recovered correctly.
            Right. The boundary $(P_1, \f t_1)$  can be recovered stable, the
            boundary $(P_2, \f t_2)$ cannot.}
    \label{fg:exlim}
    \end{psfrags}
\end{figure}

The left image in Figure \ref{fg:exlim} shows a reconstruction where data was collected for
$u \in [-20, 20]$ only (limited aperture). In this case, boundaries of some spheres
are not recovered correctly. We emphasize that this confirms theoretical predictions from
microlocal analysis \cite{LoiQui00}, namely that in order to {\em stable} recover a
boundary with tangent $\f t$ at $P_i = (u_i,v_i)$, data $G(u, \cdot)$ has to be collected
such that $( u_i-u, v_i )$ is orthogonal to $\f t$.
As illustrated in Figure \ref{fg:exlim}, this requirement is not satisfied for nearly
horizontal boundaries.

\section{Conclusion and outlook}

In this paper we presented a mathematical framework of
thermoacoustic CT with line detectors. The advantage of
using line detectors is the improved imaging resolution
compared to data acquisition based on point data. We proved that the
three dimensional imaging problem reduces to the problem of
reconstructing  the initial data of the two dimensional wave
equation from boundary measurements.

In practical applications the output of a planar detector array offers only
information from a finite aperture. Therefore, our future work will focus on
the case of limited data. Furthermore, we will investigate scanning geometries
where the array of line detectors forms different surfaces, e.g. cylinders which
allow to collect a complete data set.

\appendix

\section{The Fourier--cosine transform on the space of symmetric tempered distributions}
\label{sec:app}

In this  appendix we define the Fourier--cosine transform
on the space of symmetric tempered distributions
and prove the propositions stated in Section \ref{sec:linear}.
We will reduce the Fourier--cosine to the Fourier transform and therefore
recall the basic definitions and properties of tempered distributions and the
Fourier transform \cite{ReeSim80}.

The {\em Schwartz space} $S(\RnR)$ denotes the set of all infinitely differentiable functions
$\ph: \RnR \to \C$ that are rapidly decreasing together with all their derivatives, i.e.
\begin{equation*}
    \forall m \in \N_{0}, \fs \beta \in \N_{0}^{n+1}:
    \norm{\ph}_{m,\fs \beta} := \sup_{\f u' \in \R^{n+1}} (1 + \norm{\f u'}^{m}) \abs{D^{\fs \beta} \ph(\f u')} < \infty
\end{equation*}
The variables in $\RnR$ are denoted by $\f u' = (\f u, v)$.
We discriminate between $\f u$ and $v$ since the last component
will play a special role.
Equipped with the topology generated by the semi-norms
$\norm{\ph}_{m,\fs \beta}$, where $ m \in \N_{0} $, $\fs \beta \in \N_{0}^{n+1}$,
$S(\RnR)$ is a Fr\'{e}chet space (metrizable and locally convex).

The space $S'(\RnR)$ of {\em tempered distributions} is defined as the topological dual of  $S$, i.e.
the set of all linear continuous functionals $ T: S(\RnR) \to \C$,
equipped with the weak--${}\ast$ topology $\si(S'(\RnR), S(\RnR))$. The evaluation of a
tempered distribution is denoted by $ \ip{T}{\ph} := T(\ph)$. The mappings
$\injS: S(\RnR) \to S'(\RnR)$,
$\injL: L^1(\RnR) \to S'(\RnR)$, and
$\injC: \Czz(\RnR)  \to S'(\RnR)$,
where $\Czz(\RnR)$ denotes the set of all continuous functions $\ph$ with
$\lim_{\norm{\f u'} \rightarrow \infty} \ph(\f u')  = 0$ are well defined and injective.
In all cases a function $\ph$ is mapped to the distribution
$\inj[ \ph ]$ defined by
\[
    \ip{ \inj[\ph]}{ \psi } := \int_{\R^{n+1} } \ph(\f u') \psi(\f u') d^{n+1}\f u' \,,
\]
where $\inj$ stands either for $\injS$, $\injL$ or $\injC$.

\begin{definition}\label{defn:app}
\end{definition}
    \begin{itemize}
        \item   For $\ph : \RnI \to \C$ and $\psi : \RnR  \to \C$ define
                \begin{equation*}
                    \begin{aligned}
                        \ph_{\rm sym}: \RnR           & \to \C                                                \\
                                (\f u,v)              & \mapsto \ph_{\rm sym}(\f u, v ) := \ph(\f u, \abs{v})\,, \\
                        \psi_{-}: \RnR                & \to \C \\
                                (\f u,v)              & \mapsto \psi_{-}(\f u, v ) := \psi(\f u, -v) \,, \\
                        \psi^{-}: \RnR                & \to \C \\
                                (\f u,v)              & \mapsto \psi^{-}(\f u, v ) := \psi(-\f u, v) \,, \\
                        \psi_{-}^{-}: \RnR            & \to \C \\
                                (\f u,v)              & \mapsto \psi_{-}^{-}(\f u, v ) := \psi(-\f u, -v) \;.
                    \end{aligned}
                \end{equation*}
        \item   The {\em Fourier transform on $L^1( \RnR )$}
                is defined by
                \begin{equation*}
                    \begin{aligned}
                        & F: L^1( \RnR )               \rightarrow   \Czz(\RnR) \\
                        &       \ph                  \mapsto \ft[\ph] := \bigg(
                                                                             (\fs \la, \ka ) \mapsto
                                                                             \left( \frac{1}{2 \pi} \right)^{(n+1)/2}
                                                                             \int_{\R^n}
                                                                             \left( \int_\R \ph(\f u, v) e^{- i \ka v }  dv \right)
                                                                             e^{- i \ip{\fs \la}{\f u} }    d^{n} \f u
                                                                          \bigg)\;.
                    \end{aligned}
                \end{equation*}
        \item   The Fourier transform  $\ft: S( \RnR  ) \to S( \RnR  )$ restricted to the
                Schwartz space is well defined and a topological
                isomorphism (linear,  bijective and bicontinuous) with inverse
                $\ift [ \Phi]= \ft[\Phi]_{-}^{-}$.
        \item   The {\em Fourier transform on $S'(\RnR)$}
                \begin{equation*}
                    \begin{aligned}
                        \ftp: S'(\RnR)              & \rightarrow   S'(\RnR) \\
                            T                       & \mapsto  \ftp[T] := \big( \ph \mapsto \ip{T}{\ft[\ph]} \big)
                    \end{aligned}
                \end{equation*}
                is a topological isomorphism with $\ip{\iftp[T]}{\ph}  = \ip{T}{\ift[\ph]}$.
\end{itemize}
The definitions of the Fourier transform $\ft$ and $\ftp$ are compatible, i.e.
                for $\ph \in L^1( \RnR )$
                \begin{equation} \label{eq:FComp}
                    \ftp [\injL[\ph] ] = \injC[\ft[\ph]]
                \end{equation}
                because, by using Fubini's theorem,
                \begin{equation*}
                    \begin{aligned}
                        \ip{\ftp [\injL[\ph]]}{\psi} &= \ip{\injL[\ph]}{\ft[\psi]} \\
                                                     &= \int_{\R^{n+1}} \ph(\f x)
                                                        \left(
                                                        \left( \frac{1}{2 \pi} \right)^{(n+1)/2}
                                                        \int_{\R^{n+1}} \psi(\fs \zeta) e^{-i \ip{\f x}{\fs \zeta}} d^{n+1} \fs \zeta
                                                        \right)
                                                        d^{n+1} \f x \\
                                                     &= \int_{\R^{n+1}} \psi(\fs \zeta)
                                                        \left(
                                                        \left( \frac{1}{2 \pi} \right)^{(n+1)/2}
                                                        \int_{\R^{n+1}} \ph(\f x) e^{-i \ip{\f x}{\fs \zeta}} d^{n+1} \f x
                                                        \right)
                                                        d^{n+1} \fs \zeta \\
                                                     &= \ip{\injC [\ft[\ph]]}{\psi} \;.
                    \end{aligned}
                \end{equation*}
                In a similar way it can be shown that
                \begin{equation} \label{eq:FComp2}
                    \iftp[\injL[\ph]] = \injC[\ft[\ph]_{-}^{-}] \;.
                \end{equation}
\begin{definition}
\end{definition}
    \begin{itemize}
        \item   A tempered distribution $T \in S'(\RnR)$ is called {\em symmetric} if
                \begin{equation*}
                    \forall \ph \in S(\RnR): \ip{T}{\ph} = \ip{T}{\ph_{-}} \;.
                \end{equation*}

        \item   We define
                \begin{eqnarray*}
                    \Ssym(\RnR)               &:=& \set{\ph \in C^\infty(\RnI,\C): \ph_{\rm sym} \in S(\RnR)} \,, \\
                    \Czzsym(\RnR)             &:=& \set{\ph \in C^\infty(\RnI,\C): \ph_{\rm sym} \in \Czz(\RnR)} \,,\\
                    \Spsym(\RnR)              &:=& \set{T \in S'(\RnR,\C): T \text{ symmetric} } \;.
                \end{eqnarray*}
                The subspaces $\Ssym(\RnR)$, $\Spsym(\RnR)$ have the topologies induced
                by $S(\RnR)$ and $S'(\RnR)$.

        \item   The {\em Fourier--cosine transform $\cost$ on $L^1(\RnI)$}
                is defined by
                \begin{equation*}
                    \cost[\ph](\fs \la, \ka )
                    :=
                    2 \left( \frac{1}{2 \pi} \right)^{(n+1)/2}
                    \int_{\R^n} \left( \int_0^{\infty}
                    \ph (\f u, v) \cos(\ka v )  dv \right)  e^{- i \ip{\fs \la}{\f u} }    d^{n} \f u \;.
                \end{equation*}
        \item   The {\em Fourier--cosine transform $\costp := \ftp|_{\Spsym(\RnR)}$
                on $\Spsym(\RnR)$} denotes  the restriction of
                $\ftp: S'(\RnR) \to S'(\RnR)$ to the space $\Spsym(\RnR)$.
\end{itemize}

\begin{proposition}[$\costp$ is an isomorphism]\label{proposition:iso-s-str}
The mapping $C': \Spsym(\RnR) \to \Spsym(\RnR)$ is a
topological isomorphism and
$\ip{\icostp[T]}{\ph} = \ip{\costp[T]}{\ph^{-}}$.
\end{proposition}

\begin{proof}
The mapping
$\costp: \Spsym(\RnR) \to S'(\RnR)$ is
linear, continuous and injective. Let $T_{1},T_{2} \in \Spsym(\RnR)$
and $ \ph \in S(\RnR)$. Then
\begin{equation*}
        \begin{aligned}
                \ip{ \ftp [T_{1}] }{ \ph_{-} } &=
                \ip{ T_{1} }{ \ft[\ph_{-}] } =
                \ip{ T_{1} }{ \ft[\ph]_{-} } =
                \ip{ T_{1} }{ \ft[\ph] }=
                \ip{ \ftp[ T_{1}] }{ \ph }
                \\
                \ip{ \iftp [T_{2}] }{ \ph_{-} } &=
                \ip{ T_{2} }{ \ift[\ph_{-}] } =
                \ip{ T_{2} }{ \ift[\ph]_{-} } =
                \ip{ T_{2} }{ \ift[\ph] }=
                \ip{ \iftp[ T_{2} ] }{ \ph }
        \end{aligned}
\end{equation*}
Hence, $ \ftp [T_{1}]$ and $\iftp [T_{2}]$ are symmetric and $\costp: \Spsym(\RnR) \to \Spsym(\RnR)$
is well defined, surjective and therefore an isomorphism.
Finally,
\begin{equation*}
    \begin{aligned}
        \ip{\icostp[T_{1}]}{\ph}    & =  \ip{ \iftp[T_{1}] }{ \ph }
                                  =  \ip{ T_{1} }{ \ift[\ph] }
                                  =  \ip{ T_{1} }{ \ft[\ph_{-}^{-}] } \\
                                & =  \ip{ T_{1} }{ \ft[\ph_{-}^{-}]_{-} }
                                  =  \ip{ T_{1} }{ \ft[\ph^{-}] }
                                  =  \ip{ \ftp[T_{1}] }{ \ph^{-} }
                                  =  \ip{ \costp[T_{1}] }{ \ph^{-} }
    \end{aligned}
\end{equation*}
finishes the proof.
\end{proof}

\begin{lemma} \label{lemma:FCsym}
If $\ph \in L^1(\RnI)$, $\psi \in \Czzsym(\RnR) $, $\eta \in \Ssym(\RnR)$ then
$\injL [ \ph_{\rm sym} ], \injC [\psi_{\rm sym} ], \injS [ \eta_{\rm sym} ]
\in \Spsym(\RnR)$ and $\ft[\ph_{\rm sym}] = \cost[\ph]_{\rm sym}$.
\end{lemma}

\begin{proof}
Let $\gamma \in \set{ \ph, \psi, \eta }$ and i denote the corresponding embedding. Then
\begin{equation*}
    \begin{aligned}
        \ip{\rm i[\gamma_{\rm sym}]}{\rho}  &= \int_{\R^{n+1}} \gamma_{\rm sym}(\f u') \rho(\f u') d^{n} \f u' \\
                                            &= \int_{\R^{n}} \bigg( \int_{\R} \gamma(\f u, \abs{v}) \rho(\f u, v) dv \bigg) d^{n} \f u \\
                                            &= \int_{\R^{n}} \bigg( \int_{\R} \gamma(\f u, \abs{v}) \rho(\f u, -v) dv \bigg) d^{n} \f u
                                             = \ip{\rm i[\gamma_{\rm sym}]}{\rho_{-}} \;.
    \end{aligned}
\end{equation*}
which proves the first statement. Moreover, for $\ka \geq 0$
\begin{equation*}
    \begin{aligned}
        \cost[\ph](\fs \la, \ka)  &= 2 \left( \frac{1}{2 \pi} \right)^{(n+1)/2}
                                     \int_{\R^n} \left( \int_0^{\infty} \ph(\f u, v) \cos(\ka v) dv \right)
                                     e^{- i \ip{\fs \la}{\f u} } d^{n} \f u \\
                                  &= \left( \frac{1}{2 \pi} \right)^{(n+1)/2}
                                     \int_{\R^n} \left( \int_\R \ph_{\rm sym}(\f u, v) e ^{-i \ka v }  dv \right)
                                     e^{- i \ip{\fs \la} {\f u} } d^{n} \f u \\
                                  &= \ft[\ph_{\rm sym}](\fs \la, \ka) \;.
    \end{aligned}
\end{equation*}
$\ph_{\rm sym}$ implies that $\ft[\ph_{\rm sym}]$ is a symmetric function with respect to
the last component, i.e. $\ft[\ph_{\rm sym}] \in \Czz(\RnR)$ symmetric.
So, by extending $\cost[\ph]$ to $\cost[\ph]_{\rm sym}$ the second statement holds.
\end{proof}

\begin{proposition}[Compatibility] \label{proposition:C-comp}
If $\ph \in L^1(\RnI)$ then $\cost[\ph]_{\rm sym} \in  \Czzsym(\RnR)$. Moreover,
$\costp [ \injL[\ph_{\rm sym}] ] =  \injC [\cost[\ph]_{\rm sym}]$ and
$\icostp [ \injL[\ph_{\rm sym}] ] =  \injC [\cost[\ph]_{\rm sym}^{-}]$.
\end{proposition}

\begin{proof}
Since $\ph_{\rm sym}  \in L^1(\RnR)$ and $\ftp[\injL[\ph_{\rm sym}]] = \injC[\ft[\ph_{\rm sym}]]$ (see \eqref{eq:FComp})
the restriction of $\ftp$ to the space $\Spsym$ yields
\begin{equation} \label{eq:FCompRest1}
    \costp [ \injL[\ph_{\rm sym}] ] = \injC [ \ft [ \ph_{\rm sym} ] ] \;.
\end{equation}
Furthermore, by using \eqref{eq:FComp2}
\begin{equation} \label{eq:FCompRest2}
    \icostp[\injL[\ph]_{\rm sym}] = \injC[\ft[\ph_{\rm sym}]_{-}^{-}] = \injC[\ft[\ph_{\rm sym}]^{-}] \;.
\end{equation}
Applying Lemma \ref{lemma:FCsym} to \eqref{eq:FCompRest1} and \eqref{eq:FCompRest2} accomplishes the proof.
\end{proof}

\begin{proposition}[$\cost$ is an isomorphism] \label{proposition:iso-s}
The mapping $\cost: \Ssym(\RnR) \to \Ssym(\RnR)$
is a topological isomorphism with
\begin{equation*}
    \icost[\ph](\f u, v ) = 2 \left( \frac{1}{2 \pi} \right)^{(n+1)/2}
                            \int_{\R^n} \left( \int_\R \ph(\fs \la, \ka) \cos(\ka v )  d \ka \right)
                            e^{i \ip{\fs \la}{\f u} } d^{n} \fs \la \;.
\end{equation*}
\end{proposition}

\begin{proof}
Both statements follow  from  $\cost[\ph]_{\rm sym}  = \ft [\ph_{\rm sym}]|_{\RnI}$
and the corresponding result for $\ft$.
\end{proof}

We finally note that analogous  constructions are possible for the
space $S(\R^n \times \R^m)$ of all Schwartz functions that are
rotationally  invariant in the second component. This leads to an
extension of the Fourier-Hankel transform to a space of distributions.
In the special situation $m = n+1$ this space was denoted by
${\cal S}'_r(\R^n \times \R^{n+1})$ in \cite{And88,SchQui05}
and used to invert the spherical mean operator.

\section{Inversion of the wave equation in arbitrary dimensions.}
\label{sec:n-d}

Let $n \in \N$, $n \geq 1$, and $H \in C^\infty_0(\RnI)$.
The variables in $\RnR$ are once more denoted by $\f u' = (\f u, v)$.
We consider the problem of reconstructing the
initial data $H$ of the $(n+1)$-dimensional wave equation
\begin{equation}\label{eq:wave-np}
     \left( \frac{\partial^2}{\partial t^2} - \lap_{\f u'} \right) P = 0 \,,\quad
     (\f u, v, t)  \in \RnR  \times (0, \infty)
\end{equation}
with initial conditions
\begin{eqnarray}
  P (\f u' , 0) =  H(\f u, \abs{v})\,,  & \f u' \in \RnR \,, \label{eq:ini-np1} \\
  \frac{\partial P}{\partial t} (\f u' , 0 ) = 0  \,, & \f u' \in \RnR  \label{eq:ini-np2}
\end{eqnarray}
from data $G(\f u , t) := P (\f u , v=0, t)$.
Here $\lap_{\f u'}$ denotes the Laplacian in $\RnR$.

\begin{theorem}\label{thm:fif-n}
Let $H \in C^\infty_0 (\R^{n} \times (0, \infty))$, $P$ the solution
of \req{wave-np}--\req{ini-np2} and $G := P|_{\R^n \times \set{0}
\times (0, \infty)}$ denote the restriction of $P$ to $\R^n \times
\set{0} \times (0, \infty)$.

Then  there exists $\bar{G} \in L^1(\R^n \times (0, \infty))$, such
that  $\costp[\injC[G_{\rm sym}]] = \injL[ \bar{G}_{\rm sym}]$.
Moreover,
  \begin{equation*}
    \bar{G}  (\fs \la, \om ) =   \lim_{K \to \infty}
                                        2 \left( \frac{1}{2\pi} \right)^{(n+1)/2}
                                        \int_0^K
                                        \left(
                                                \int_{\R^n} G ( \f u, t ) e^{ -i \ip{\fs \la}{\f u}}  d^n \f u
                                        \right)
                                        \cos( \om t )  \, dt
\end{equation*}
exists point--wise for all $\omega > \lambda$ and
\begin{equation*}
    \cost [ H ] ( \fs \la, \ka  ) =
     \bar{G} \left(  \fs \la, \sqrt{ \ka^2 + \norm{\fs \la}^2 } \right)
     \cdot \frac{\ka}{\sqrt{ \ka^2 + \norm{\fs \la}^2  }} \,, \quad \ka > 0 \;.
\end{equation*}
\end{theorem}

The proof of Theorem \ref{thm:fif-n} is analogous to the corresponding
result in  $n+1 = 2$ dimensions, as presented in Section \ref{sec:linear}.

\section*{Acknowledgement}
This work has been supported by the Austrian Science Fund (FWF),
Project Y123--INF and project P18172--N02.

\bibliography{TCTline}

\begin{thebibliography}{10}

\bibitem{AgrQui96}
M.~Agranovsky and E.~Quinto.
\newblock Injectivity sets for the {R}adon transform over circles and complete
  systems of radial functions.
\newblock {\em Journal of Functional Analysis}, 139(2):383--414, 1996.

\bibitem{art05:AmbKuc}
G.~Ambartsoumian and P.~Kuchment.
\newblock On the injectivity of the circular radon transform.
\newblock {\em Inverse Problems}, 21(2):473--485, 2005.

\bibitem{And88}
L.-E. Andersson.
\newblock On the determination of a function from spherical averages.
\newblock {\em SIAM Journal on Mathematical Analysis}, 19(1):214--232, 1988.

\bibitem{AKO01}
V.G. Andreev, A.A. Karabutov, and A.A. Oraevsky.
\newblock Detection of ultrawide-band ultrasound pulses in optoacoustic
  tomography.
\newblock {\em IEEE Transactions on Ultrasonics, Ferroelectrics, and Frequency
  Control}, 50:1383--1390, 2003.

\bibitem{BurEtAl05}
P.~Burgholzer, C.~Hofer, G.~Paltauf, M.~Haltmeier, and O.~Scherzer.
\newblock Thermoacoustic tomography with integrating area and line detectors.
\newblock {\em IEEE Transactions on Ultrasonics, Ferroelectrics, and Frequency
  Control}, 52:1577--1583, 2005.

\bibitem{CoxEtAl04}
B.~T. Cox, E.~Z. Zhang, Laufer~J. G., and Beard~P. C.
\newblock Fabry perot polymer film fibre-optic hydrophones and arrays for
  ultrasound field characterisation.
\newblock {\em Journal of Physics: Conference Series}, 1:32--37, 2004.

\bibitem{FinPatRak04}
D.~Finch, S.~Patch, and Rakesh.
\newblock Determining a function from its mean values over a family of spheres.
\newblock {\em SIAM Journal on Mathematical Analysis}, 35(5):1213--1240, 2004.

\bibitem{GusKar96}
V.E. Gusev and A.A. Karabutov.
\newblock {\em Laser Optoacoustics}.
\newblock Institute of physics, New York, 1993.

\bibitem{art04:HalEtAl}
M.~Haltmeier, O.~Scherzer, P.~Burgholzer, and G.~Paltauf.
\newblock Thermoacoustic computed tomography with large planar receivers.
\newblock {\em Inverse Problems}, 20:1663--1673, 2004.

\bibitem{HSS05}
M.~Haltmeier, T.~Schuster, and O.~Scherzer.
\newblock Filtered backprojection for thermoacoustic computed tomography in
  spherical geometry.
\newblock {\em Mathematical Methods in the Applied Sciences},
  28(16):1919--1937, 2005.

\bibitem{Jen99}
J.~A. Jensen.
\newblock A new calculation procedure for spatial impulse responses in
  ultrasound.
\newblock {\em Journal of the Acoustical Society of America}, 105:3266--3274,
  1999.

\bibitem{Joh82}
F.~John.
\newblock {\em Partial Differential Equations}.
\newblock Springer, Berlin-New York, 1982.

\bibitem{KoeBea03}
K.~P. K\"ostli and P.~C. Beard.
\newblock Two-dimensional photoacoustic imaging by use of fourier-transform
  image reconstruction and a detector with an anisotropic response.
\newblock {\em Applied Optics}, 42(10), 2003.

\bibitem{KoeEtAl01}
K.~P. K\"ostli, D.~Frauchinger, J.~J. Niederhauser, G.~Paltauf, H.~W. Weber,
  and M.~Frenz.
\newblock Optoacoustic imaging using a three-dimensional reconstruction
  algorithm.
\newblock {\em IEEE Journal on selected topics in Quantum electronics},
  7(6):918--923, 2001.

\bibitem{KruEtAl03}
Robert~A. Kruger, Jr. William L.~Kiser, Daniel~R. Reinecke, and Gabe~A. Kruger.
\newblock Thermoacoustic computed tomography using a conventional linear
  transducer array.
\newblock {\em Medical Physics}, 30(5):856--860, 2003.

\bibitem{LoiQui00}
A.K. Louis and E.T. Quinto.
\newblock Local tomographic methods in sonar.
\newblock In {\em Surveys on solution methods for inverse problems}, pages
  147--154. Springer, Vienna, 2000.

\bibitem{Nat86}
F.~Natterer.
\newblock {\em The Mathematics of Computerized Tomography}.
\newblock Wiley, Chichester, 1986.

\bibitem{Pal04}
Victor Palamodov.
\newblock {\em Reconstructive integral geometry}, volume~98 of {\em Monographs
  in Mathematics}.
\newblock Birkh\"auser Verlag, Basel, 2004.

\bibitem{PNHB07}
G.~Paltauf, R.~Nuster, M.~Haltmeier, and P.~Burgholzer.
\newblock Thermoacoustic computed tomography using a mach-zehnder
  interferometer as acoustic line detector.
\newblock 46:3352--3358, 2007.

\bibitem{PalSchGus96}
G.~Paltauf, H.~Schmidt-Kloiber, and H.~Guss.
\newblock Light distribution measurements in absorbing materials by optical
  detection of laser-induced stress waves.
\newblock {\em Applied Physics Letters}, 69:1526--1528, September 1996.

\bibitem{ReeSim80}
M.~Reed and B.~Simon.
\newblock {\em Functional analysis}, volume~I of {\em Methods of modern
  mathematical physics}.
\newblock Academic Press, New York-London, 1980.

\bibitem{SchQui05}
Thomas Schuster and Eric~Todd Quinto.
\newblock On a regularization scheme for linear operators in distribution
  spaces with an application to the spherical {R}adon transform.
\newblock {\em SIAM Journal on Applied Mathematics}, 65(4):1369--1387, 2005.

\bibitem{Tri01}
Khalifa Trim{\`e}che.
\newblock {\em Generalized harmonic analysis and wavelet packets}.
\newblock Gordon and Breach Science Publishers, Amsterdam, 2001.

\bibitem{WPKXSW03}
X.D. Wang, G.~Pang, Y.J.~Ku, X.Y. Xie, G.~Stoica, and L.-H.V. Wang.
\newblock Noninvasive laser-induced photoacoustic tomography for structural and
  functional in vivo imaging of the brain.
\newblock {\em Nature Biotechnology}, 21:803--806, 2003.

\bibitem{XuMFenWan03}
M.~Xu, D.~Feng, and L.-H. Wang.
\newblock Time-domain reconstruction algorithms and numerical simulations for
  thermoacoustic tomography in various geometries.
\newblock {\em IEEE Transactions on Biomedical Engineering}, 50:1086-- 1099,
  2003.

\end{thebibliography}
\bibliographystyle{plain}
\end{document}